\documentclass[12pt,reqno]{amsart}
\usepackage{amsthm,amsmath,amsfonts,amssymb}
\usepackage{hyperref,pstricks}
\usepackage{color}
\usepackage{amsmath,amsfonts}
\usepackage{mathrsfs}
\usepackage{graphicx}
\usepackage{pstricks-add,amsmath}

\newcommand{\R}{\mathbb{R}}
\newcommand{\N}{\mathbb{N}}
\newtheorem{theorem}{Theorem}

\newtheorem{conjecture}[theorem]{Conjecture}

\newtheorem{definition}[theorem]{Definition}

\newtheorem{lemma}[theorem]{Lemma}

\newtheorem{proposition}[theorem]{Proposition}
\newtheorem{remark}[theorem]{Remark}

\newcommand{\cqd}{\hfill \rule{5pt}{5pt}}
\newtheorem{theorem*}{}

\begin{document}

\title[]
{A non-autonomous scalar one-dimensional dissipative parabolic problem: The description of the dynamics}

\author[R. C. D. S. Broche]{Rita de C\'assia D. S. Broche$^*$}\thanks{$^*$Work partially done while the author visited the the University of S\~{a}o Paulo at S\~{a}o Carlos - SP, the author wishes to thank the Universidade Federal de Lavras for all the support concerning that visit.}
\address[R. C. D. S. Broche]{Departamento de Ci\^encias Exatas, Universidade Federal de Lavras, Caixa Postal 3037, 37200-000 Lavras MG, Brazil}
\email{ritabroche@dex.ufla.br}

\author[A. N. Carvalho]{Alexandre N. Carvalho$^\dag$}\thanks{$^\dag$Research partially supported by CNPq 303929/2015-4, CAPES/DGU 267/2008 and FAPESP 2003/10042-0, Brazil}

\address[A. N. Carvalho]{Departamento de
Matem\'atica, Instituto de Ci\^encias Mate\-m\'a\-ti\-cas e de
Computa\c{c}\~{a}o, Universidade de S\~{a}o Paulo-Campus de S\~{a}o
Carlos, Caixa Postal 668, 13560-970 S\~{a}o Carlos SP, Brazil}
\email{andcarva@icmc.usp.br}

\author[J. Valero]{Jos\'{e} Valero$^\ddag$}\thanks{$^\ddag$ Research partially supported by Spanish Ministry of Economy and Competitiveness and FEDER, projects MTM2015-63723-P and MTM2016-74921-P, and
by Junta de Andaluc\'{\i}a (Spain), project P12-FQM-1492}

\address[J. Valero]{Centro de Investigaci\'{o}n Operativa, Universidad Miguel Hern\'{a}ndez de Elche, Avda. Universidad S/n, 03540-Elche (Alicante), Spain}

\email{jvalero@umh.es}

\subjclass[2010]{35K91 (35B32\, 35B41\, 37L30) }
\keywords{Non-autonomous Chafee-Infante Problem, Non-autonomous dynamical system, Parabolic equations, Pullback attractors, Uniform Attractors, Gradient structure}

\begin{abstract}
The purpose of this paper is to give a characterization of the structure of non-autonomous attractors of the problem
$u_t= u_{xx} + \lambda u - \beta(t)u^3$ when the parameter $\lambda > 0$ varies. Also, we answer a question proposed in \cite{Alexandre}, concerning the complete description of the structure of the pullback attractor of the problem when $1<\lambda <4$ and, more generally, for $\lambda \neq N^2$, $2\leq N\in \N$. We construct global bounded solutions , ``non-autonomous equilibria", connections between the trivial solution these ``non-autonomous equilibria" and characterize the $\alpha$-limit and $\omega$-limit set of global bounded solutions. As a consequence, we show that the global attractor of the associated skew-product flow has a gradient structure. The structure of the related pullback an uniform attractors are derived from that.

\end{abstract}

\maketitle

\section{Introduction}

Consider the semilinear parabolic problem
\begin{equation}\label{naci}
\begin{split}
&u_t=u_{xx} + \lambda u - \beta(t) u^3, \quad t>s, \ x\in (0,\pi),\\
&u(0,t)=u(\pi,t)=0, \quad t\geq s,\\
&u(x,s)=u_0(x), \quad u_0\in H^1_0(0,\pi),
\end{split}
\end{equation}
where $\lambda > 0$ and $\beta:\R\to \R$ is a continuously differentiable function with $0<\beta_1\leq\beta(t)\leq \beta_2$, for all $t\in \R$ and some real constants $\beta_1$ and $\beta_2$.

It is well known (see \cite{Henry}) that, for each $u_0\in H^1_0(0,\pi)$ and $s\in \R$, there is a unique $u(\cdot, s, u_0)\in C([s,\infty),H^1_0(0,\pi))$ which is a mild solution for \eqref{naci}.
This solution is shown to be classical for each $t>s$ and if $\mathcal{P}=\{(t,s)\in \R^2: t\geq s\}$ the map
$$
\mathcal{P}\times H^1_0(0,\pi) \ni ((t,s),u_0) \mapsto u(t,s,u_0)\in H^1_0(0,\pi)
$$
is continuous.

\medskip

With the above notation, for each $(t,s)\in \mathcal{P}$ and $u_0\in H^1_0(0,\pi)$, define $
T_\beta(t,s)u_0= u(t,s,u_0)$. It is clear that
\begin{enumerate}
\item $T_\beta(t,t)=I$
\item $T_\beta(t,\tau)T_\beta(\tau,s)=T_\beta(t,s)$, for all $s \leq \tau \leq t$ and
\item $\mathcal{P}\times H^1_0(0,\pi) \ni (t,s,u_0)\mapsto T_\beta(t,s)u_0\in H^1_0(0,\pi)$ is continuous.
\end{enumerate}

\medskip

A family of operators $\{T_\beta(t,s):(t,s)\in \mathcal{P}\}$ with the above properties is called an \emph{evolution process} in $H^1_0(0,\pi)$.

We are interested in the description of the asymptotic dynamics of the solutions of \eqref{naci}, in particular, we are interested in the family of asymptotic sets called pullback attractors. Next we introduce the basic notions needed to define pullback attractors starting with the notions of invariance and pullback attraction (see \cite{CLR12}):

\medskip

\emph{A family $\{A(t):t\in \R\}$ of subsets of $H^1_0(0,\pi)$ is {\bf invariant} under the action of the evolution process $\{T_\beta(t,s): (t,s)\in \mathcal{P}\}$ if $T_\beta(t,s)A(s)=A(t)$ for all $(t,s)\in \mathcal{P}$}.

\medskip

Recall that a continuous function $\xi:\R\to H^1_0(0,\pi)$ is a global solution for \eqref{naci} or, equivalently, for the evolution process $\{T_\beta(t,s):(t,s)\in \mathcal{P}\}$, if $T_\beta(t,s)\xi(s)=\xi(t)$ for all $(t,s)\in \mathcal{P}$.

\medskip

\emph{Given $B_0,B\subset H^1_0(0,\pi)$, we say that $B_0$ {\bf pullback-attracts $B$ at time $t$} under the action of the evolution process $\{T_\beta(t,s): (t,s)\in \mathcal{P}\}$ if
$$
\lim_{s\to-\infty} {\rm dist}_H(T_\beta(t,s)B,B_0)=0,
$$}
where ${\rm dist}_H$ denotes the Hausdorff semidistance in $H^1_0(0,\pi)$. We are now ready to define pullback attractors

\begin{definition}
We say that a family $\{\mathscr{A}_\beta(t):t\in \R\}$ is a {\bf pullback-attractor} for $\{T_\beta(t,s):(t,s)\in \mathcal{P}\}$ if it is invariant, $\mathscr{A}_\beta(t)$ is compact,  $\mathscr{A}_\beta(t)$ pullback-attracts bounded subsets of $H^1_0(0,\pi)$ at time $t$ for each $t\in \R$ and $\{\mathscr{A}_\beta(t):t\in \R\}$ is the minimal closed family with this pullback-attracting property; that is, each family of closed sets $\{C(t):t\in \R\}$ such that $C(t)$ pullback-attracts bounded subsets of $H^1_0(0,\pi)$ at time $t$, for each $t\in \R$, must satisfy $\mathscr{A}_\beta(t)\subset C(t)$ for each $t\in \R$.

If $\displaystyle{\cup_{t\leq 0}\mathscr{A}_\beta(t)}$ is bounded, it is easy to see that the minimality condition is automatically satisfied.

\end{definition}

Under the requirement that $\displaystyle{\cup_{t\leq 0}\mathscr{A}_\beta(t)}$ be bounded, the pullback attractor has the following characterization
$$
\mathscr{A}_\beta(t)=\left\{\xi(t): \xi:\R\to H^1_0(0,\pi) \hbox{ is a backwards bounded global solution of } \eqref{naci}\right\}.
$$

It is not difficult to prove that $\{T_\beta(t,s):(t,s)\in \mathcal{P}\}$ has a pullback attractor $\{\mathscr{A}_\beta(t):t\in \R\}$ (see \cite{LS}) with the property that  $\displaystyle{\cup_{t\leq 0}\mathscr{A}_\beta(t)}$ is bounded for each $t\in \R$.

When $\beta(t)\equiv \beta=\,$const we have that $T_\beta(t,s)=T_\beta(t-s,0)$ (the evolution depends only on the elapsed time) and the evolution processes is said to be autonomous. In this case
$$
S_\beta(t)=T_\beta(t,0), \ t\geq 0
$$
is a semigroup; that is,

\begin{itemize}
\item[(i)] $S_\beta(0)=I$,
\item[(ii)] $S_\beta(t+s)=S_\beta(t)S_\beta(s)$, for all $t,s\geq 0$, and
\item[(iii)] $\R^+\times H^1_0(0,\pi)\ni (t,u_0)\mapsto S_\beta(t)u_0\in H^1_0(0,\pi)$ is continuous.
\end{itemize}

For $\beta=$const, the autonomous evolution process $\{T_\beta(t,s): t\geq s\}$ has a pullback attractor $\{\mathscr{A}_\beta(t):t\in \R\}$ if and only if $\mathscr{A}_\beta(t)=\mathscr{A}_\beta$ for all $t\in \R$, and the associated semigroup $\{S_\beta(t):t\geq 0\}$ has a global attractor $\mathscr{A}_\beta$; that is,
$$
S_\beta(t)\mathscr{A}_\beta=\mathscr{A}_\beta \hbox{ for all }t\geq 0 \quad \hbox{ and } \quad \lim_{t\to\infty} {\rm dist}_H(S_\beta(t)B,\mathscr{A}_\beta)=0
$$
for all $B\subset H^1_0(0,\pi)$ bounded. A compact subset $\mathscr{A}_\beta$ of $H^1_0(0,\pi)$ with the above properties is called a global attractor for the semigroup $\{S_\beta(t):t\geq 0\}$ in $H^1_0(0,\pi)$. It is easy to see that if $\xi:\R\to H^1_0(0,\pi)$ is a global bounded solution for $\{S_\beta(t):t\geq 0\}$ then $\overline{\xi(\R)}\subset \mathscr{A}_\beta$.

The aim of this paper is to reveal the little we know about the internal dynamics of the pullback attractor $\{\mathscr{A}_\beta(t):t\in \R\}$ for $\beta(\cdot)$ not necessarily close to a constant.

\medskip

Very little is known about the internal dynamics of pullback attractors for non-autonomous evolution processes, even for very simple models. The aim of this work is to reveal a little of the dynamics of such processes. The choice of \eqref{naci} for this study is related to the fact that, in the autonomous case, it is the infinite-dimensional model for which the asymptotic dynamics is best understood. We will take advantage of the work \cite{Alexandre}, where the existence of the non-autonomous equilibria has been established, and will establish some of the connections between these non-autonomous equilibria, inspired by the description given for the autonomous case, which we state next.

\medskip


Consider the classical autonomous Chafee-Infante problem (see \cite{CI,Hale,Henry})
\begin{equation}\label{aci}
\begin{split}
&u_t=u_{xx} + \lambda u - \beta u^3, \quad t>0, \ x\in (0,\pi),\\
&u(0,t)=u(\pi,t)=0, \quad t\geq 0\\
&u(x,0)=u_0(x), \quad u_0\in H^1_0(0,\pi),
\end{split}
\end{equation}
$\lambda\in (0,\infty)$ and $\beta>0$.

First note that the semigroup $\{S_\beta(t): t\geq 0\}$ associated with \eqref{aci} is gradient; that is, there is a continuous function $V:H^1_0(0,\pi) \to \R$ such that $[0,\infty)\ni t \mapsto V(S_\beta(t)\phi)\in \R$ is non-increasing and if $V(S_\beta(t)\phi)=V(\phi)$ for all $t\geq 0$ we must have that $\phi$ is an equilibrium solution for \eqref{aci}; that is, an element of the set $\mathcal{E}$ of solutions of the boundary value problem
\begin{equation}\label{eq-ci}
\begin{split}
&\phi'' + \lambda \phi - \beta \phi^3=0, \ x\in (0,\pi),\\
&\phi(0)=\phi(\pi)=0.
\end{split}
\end{equation}
A function $V:H^1_0(0,\pi)\to \R$ with these properties is called a Lyapunov function for \eqref{aci} in $H^1_0(0,\pi)$.
In fact, if $V:H^1_0(0,\pi)\to \R$ is defined by
$$
V(\phi)=\frac{1}{2} \int_0^\pi \phi'(x)^2 dx - \frac{\lambda}{2}\int_0^\pi  \phi(x)^2 dx + \frac{\beta}{4} \int_0^\pi \phi(x)^4 dx
$$
then, for $u_0\in H^1_0(0,\pi)$ and $u(t,x)=S_\beta(t)u_0(x)$,
$$
\frac{d}{dt} V(S_\beta(t)u_0) = -\int_0^\pi u_t(t,x)^2 dx.
$$
Hence, $V$ is a Lyapunov function for \eqref{aci}.

The most elementary bounded solutions of $\{S_\beta(t):t\geq 0\}$ are the equilibria and $\mathcal{E}\subset \mathscr{A}_\beta$. Next consider the global bounded solutions that converge to an equilibria as $t\to -\infty$. Given $\phi\in \mathcal{E}$ we define the unstable set of $\phi$ as $
W^u(\phi)=\{ u\in H^1_0(0,\pi): \hbox{ there exists global solution }$ $\xi:\R \to H^1_0(0,\pi) \hbox{ of } \{S_\beta(t):t\geq 0\} \hbox{ such that } \xi(0)=u \hbox{ and }\lim_{t\to -\infty} \xi(t)=\phi\}$. It is clear that $W^u(\phi)\subset \mathscr{A}_\beta$ for all $\phi\in \mathcal{E}$.

It is well known that a gradient semigroup $\{S_\beta(t):t\geq 0\}$ with a global attractor $\mathscr{A}_\beta$ and such that the $\mathcal{E}$ has finitely many elements (we will see that this is the case for $\{S_\beta(t):t\geq 0\}$) has the property that, given a global bounded solution $\xi:\R\to H^1_0(0,\pi)$ for $\{S_\beta(t):t\geq 0\}$, there are $\phi^-_\xi,\phi^+_\xi \in \mathcal{E}$ such that $ \phi^-_\xi \stackrel{t\to -\infty}{\longleftarrow}\xi(t)\stackrel{t\to +\infty}{\longrightarrow} \phi^+_\xi$ and $V(\phi^-_\xi)>V(\phi^+_\xi)$. This immediately implies that
$\mathscr{A}_\beta=\bigcup_{\phi\in \mathcal{E}} W^u(\phi)$.

It is proved in \cite{CI} that, if $0<\lambda \leq 1$ the only possible solution of \eqref{eq-ci} is $\phi_{0}\equiv 0$, if $\lambda\in (1,4]$, there will be exactly three elements $\phi_{0,\beta}$, $\phi_{1,\beta}^+$ (positive in $(0,\pi)$) and $\phi_{1,\beta}^-$ (negative in $(0,\pi)$) in $\mathcal{E}$. For $\lambda\in (4,9]$ we will have two additional solutions $\phi_{3,\beta}^\pm$ which change sign at $x=\pi/2$ and this procedure yields a sequence of pitchfork bifurcations at $\lambda_n=n^2$ (see \cite{CI,CLR12} for details). If $\lambda\in (n^2,(n+1)^2]$ the set $\mathcal{E}$ of solutions of \eqref{eq-ci} has exactly $2n+1$ elements
$$
\phi_{0},\phi_{j,\beta}^\pm, \quad 1\leq j \leq n.
$$
The solutions $\phi_{j,\beta}^\pm$ of \eqref{eq-ci} bifurcate from $\phi_0\equiv 0$ at $\lambda=j^2$. $\phi_{j,\beta}^\pm$ has $j+1$ zeroes in $[0,\pi]$ ($x_i=\frac{i\pi}{j}$, $0\leq i\leq j$). As a consequence of that, the set $\mathcal{E}$ will always be finite and the attractor $\mathscr{A}_\beta$ will always be $\mathscr{A}_\beta=\bigcup_{\phi\in \mathcal{E}} W^u(\phi)$. Furthermore, for any $u_0\in H^1_0(0,\pi)$ $S_\beta(t)u_0 \stackrel{t\to\infty}{\longrightarrow} \phi$ for some $\phi\in \mathcal{E}$.

\begin{remark}
In the non-autonomous case, this analysis is no longer available and we must seek different ways to find the solutions that should play the role of equilibria. These solutions indeed exist but they are obtained in a completely different manner $($see \cite{CLR12-1}$)$.
\end{remark}

When $\lambda\neq j^2$, for all $j\in \N$, all equilibria in $\mathcal{E}$ are hyperbolic (see \cite{CI}). In this case making a small ($C^1$ small) autonomous perturbation of the righthand side (even if the perturbation contains gradient terms) of \eqref{aci} will result (see \cite{CL09}) a perturbed problem which will have a global attractor $\tilde{\mathscr{A}}_{\beta}$ with the same number of equilibria ($\tilde{\mathcal{E}}$ is the new set of equilibria) all of them hyperbolic, and ${\tilde{\mathscr{A}}}_\beta=\bigcup_{\phi\in \tilde{\mathcal{E}}} W^u(\tilde{\phi})$. It was proven in \cite{ACCL} that the perturbed system will also have a Lyapunov function.

The solutions $u: [0,\pi]\times [0,\infty) \to \R$ of \eqref{naci} or a linear version of it have another striking property called ``Lap Number''  (see \cite{Angenent,Matano}). This property is roughly described as follows: if $t_1>t_2>0$, ``the number of times $u(\cdot,t_1)$ vanishes  in the interval $[0,\pi]$ is at most the number of times $u(t_2,\cdot)$ vanishes in the interval $[0,\pi]$''. As a consequence of this nice property (which we will properly state later in the paper) we have that:

\begin{itemize}

\newcommand{\ti}{\;\;\makebox[0pt]{$\top$}\makebox[0pt]{$\cap$}\;\;}

\item If $\xi:\R\to H^1_0(0,\pi)$ is a global solution for $\{S_\beta(t):t\geq 0\}$ such that $\xi(t)\stackrel{t\to \pm\infty}{\longrightarrow} \phi^\pm$, then
$W^u(\phi^-)\ti W^s(\phi^+)$ and the semiflow $\{S_\beta(t):t\geq 0\}$ defined by \eqref{aci} is Morse-Smale (see \cite{Henry2,Angenent2}).

\end{itemize}

\begin{itemize}

\item The connections between equilibria are given by the following diagram (see \cite{FRJDE1996})
\end{itemize}

\begin{figure}
\begin{pspicture}(-3.5,1)(12,4)
\psset{xunit=.55cm,yunit=.6cm}

\uput[180](0.1,4){\tiny{$\phi_0$}}
\uput[180](0.48,3.98){\tiny{$\bullet$}}
\psline{->}(0,4)(2,6)
\uput[180](2.7,6.7){\tiny{$\phi_n^+$}}
\psline{->}(0,4)(2,2)
\uput[180](2.7,1.5){\tiny{$\phi_n^-$}}
\psline{->}(2,2)(6,6)
\uput[180](6.7,6.7){\tiny{$\phi_{n-1}^+$}}
\psline{->}(2,6)(6,2)
\uput[180](6.7,1.5){\tiny{$\phi_{n-1}^-$}}
\psline{->}(2,6)(6,6)
\psline{->}(2,2)(6,2)
\psline[linestyle=dashed,dash=3pt 2pt]{->}(6,6)(7,5)
\psline[linestyle=dashed,dash=3pt 2pt]{->}(6,2)(7,3)
\psline[linestyle=dashed,dash=3pt 2pt]{->}(6,6)(7,6)
\psline[linestyle=dashed,dash=3pt 2pt]{->}(6,2)(7,2)
\psline[linestyle=dashed,dash=3pt 2pt]{->}(9,2)(10,2)
\psline[linestyle=dashed,dash=3pt 2pt]{->}(9,6)(10,6)
\psline[linestyle=dashed,dash=3pt 2pt]{->}(9,3)(10,2)
\psline[linestyle=dashed,dash=3pt 2pt]{->}(9,5)(10,6)
\psline{->}(10,2)(14,6)
\uput[180](10.7,6.7){\tiny{$\phi_{2}^+$}}
\psline{->}(10,6)(14,2)
\uput[180](10.7,1.5){\tiny{$\phi_{2}^-$}}
\psline{->}(10,6)(14,6)
\uput[180](14.7,6.7){\tiny{$\phi_{1}^+$}}
\psline{->}(10,2)(14,2)
\uput[180](14.7,1.5){\tiny{$\phi_{1}^-$}}

\end{pspicture}
\caption{Diagram of connections}
\end{figure}
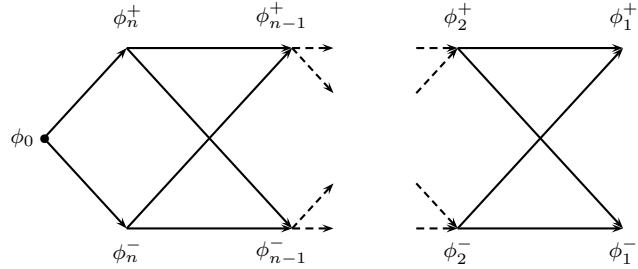
\bigskip
\bigskip
\bigskip

The above diagram has to be interpreted in the following way. If there is a sequence of oriented
segments connecting $\phi$ to $\psi$ then there exists a global solution $\xi_{\phi,\psi}:\R\to H^1_0(0,\pi)$ such that $
\phi\stackrel{t\to-\infty}{\longleftarrow} \xi_{\phi,\psi}(t)\stackrel{t\to\infty}{\longrightarrow} \psi$. Consequently, there are connections from $\phi_0$ to any other equilibria and no connection from $\phi_j^\pm$, $1\leq j\leq n$, to $\phi_0$. There are connections from $\phi_n^+$ ($\phi_n^-$) to all equilibria except to $\phi_0$ and $\phi_n^-$ ($\phi_n^+$) and so on. As a general rule, there are connections between one equilibrium $\phi$ and another equilibrium $\psi$ if $\phi$ vanishes ``more times'' than $\psi$ in $[0,\pi]$.

\begin{remark} At this point we remark that, as a consequence of the structure of attractors and of the ``Lap Number Property'', the only global bounded solutions of \eqref{aci} not lying in the unstable manifold of zero for which the zeroes in $[0,\pi]$ do not move as $t$ varies are the equilibrium solutions (the elements of $\mathcal{E}$). This will be the key to find the solutions that will play the role of equilibria when $\beta(t)$ is not close to a constant $($see \cite{CLR12-1}$)$.
\end{remark}

As a consequence of the transversality, the perturbed attractors $\tilde{\mathscr{A}}_\beta$ will have exactly the same structure as ${\mathscr{A}}_\beta$, that is, pictorially ``the connections between equilibria are the same'' in the perturbed or in the unperturbed attractor.

\medskip

At this point it is natural to ask what happens in the non-autonomous case of \eqref{naci}.
Of course, with so much structure for the global attractor of \eqref{aci} and when $\lambda\neq j^2$, for all $j\in \N$, we intuitively guess that much of the ``structure'' of the ``attractors'' must remain the same. This is, in fact, the case as a consequence of the results in \cite{CL09,ACCL,BCLR} when $\beta(t)$ is a small non-autonomous perturbation of a constant. In this paper we will prove that much of this structure remains the same even when $\beta(t)$ is not a small perturbation of a constant.

\begin{definition}\label{s-pm-beta}
$\mathcal{S}^\pm(\beta)$ denotes the set of all functions $\gamma:\R\to \R$ obtained as, uniform in bounded sets, limits of sequences $\beta(\cdot + t_n)$ with $t_n\stackrel{n\to\infty}{\longrightarrow}\pm\infty$, respectively.
\end{definition}

In fact from \cite{Alexandre}, for $N^2<\lambda< (N+1)^2$, there are $2N$ ``non-autonomous" equilibria $\xi_{j,\beta}^{\pm}$, $1\leq j\leq N$, where the index $\beta$ indicates the dependence on $\beta$.

We prove that for any $u_0\in H^1_0(0,\pi)$ the solution converges to $\{\xi_{j,\gamma}^{\pm}(t):\gamma \in \mathcal{S}^+(\beta),\; t \in \mathbb{R}\}$ as $t\to \infty$. We also prove that if there is a global bounded solution through $u_0$, it converges to $\{\xi_{j,\gamma}^{\pm}(t):\gamma \in \mathcal{S}^-(\beta),\; t \in \mathbb{R}\}$ as $t\to -\infty$. We also prove that there are solutions connecting the zero equilibrium solution to all non-autonomous equilibria.

In Section \ref{basic} we recall some basic properties of solutions of \eqref{naci} and the construction (see \cite{RB-VL}) of a global non-degenerate solution in the positive cone for $\lambda>1$ (the first non-autonomous equilibria). Section \ref{S2} is dedicated to the characterisation of the $\omega-$ and $\alpha-$limit of solutions inspired by \cite{chen}. In Section \ref{S3} we prove that, for $1<\lambda<4$ all global bounded solutions are forwards asymptotic to one of the non-autonomous equilibria and backwards asymptotic to the zero equilibrium solution. In Section \ref{valero} we characterize the $\omega-$limit and $\alpha-$limit of solutions of \eqref{naci} giving a characterization of the pullback attractor. In Section \ref{S4} we show that the zero equilibrium solution connects with all non-autonomous equilibria. Finally, in Section \ref{comments} we make some general comments about the gradient structure of the corresponding skew-product attractors and state a conjecture on our beliefs with respect to some additional structure of the dynamics of \eqref{naci}.

\section{Basic facts}\label{basic}

In this section we collect some basic facts and known results which will be used throughout.

\subsection{Special properties of solutions of parabolic problems}

In this section we collect some special properties of solutions of scalar one-dimensional parabolic equations.

\begin{description}

\item[(a) Scaling] If $\beta_1$ and $\beta_2$ are positive numbers and $u_{\beta_1}$ is a solution of
\begin{equation}\label{aci-beta1}
\begin{split}
&u_t=u_{xx} + \lambda u - \beta_1 u^3, \quad t>0, \ x\in (0,\pi),\\
&u(0,t)=u(\pi,t)=0, \quad t\geq 0,
\end{split}
\end{equation}
then $u_{\beta_2}=\left(\frac{\beta_1}{\beta_2}\right)^\frac12 u_{\beta_1}$ is a solution of
\begin{equation}\label{aci-beta2}
\begin{split}
&u_t=u_{xx} + \lambda u - \beta_2 u^3, \quad t>0, \ x\in (0,\pi),\\
&u(0,t)=u(\pi,t)=0, \quad t\geq 0.
\end{split}
\end{equation}

\item[(b) Symmetry]

If $u_0(x) = \pm u_0(\pi-x)$ and $u(t,s,u_0)(x):=T_\beta(t,s)u_0(x)$, then
$$
u(t,s,u_0)(x)=\pm u(t,s,u_0)(\pi-x),\ \hbox{ for all } t\geq s, \ x\in [0,\pi],
$$
and the zeroes of $u(t,s,u_0)(\cdot)$ are symmetric with respect to $\pi/2$.
In particular, if $u_0(x)=-u_0(\pi-x)$, then $\pi/2$ is a zero for $u_0(\cdot)$ and $u(t,s,u_0)(\pi/2)=0$ for all $t\geq 0$.

\item[(c) Comparison]

If $u_1\geq u_2$, then $T_\beta(t,s)u_1\geq T_\beta(t,s)u_2$ for all $t\geq s$ and if $u_0\geq 0$, since $\beta_1 \leq \beta(t) \leq \beta_2$,
$$
S_{\beta_2}(t-s)u_0\leq T_\beta(t,s)u_0 \leq S_{\beta_1}(t-s)u_0
$$
for all $t\geq s$.

\item[(d) Lap Number]


Now, let $\mathcal{C}_P=\{u\in \mathcal{C}(\R):u \hbox{ is $2\pi$ periodic}\}$ and define the map $\ell: C_P \to \mathbb{N} \cup \{\infty\}$, by
\begin{center}
$\ell(w)=$ the number of points in $[-\pi,\pi]$ for which $w(x)=0$.
\end{center}
The following result is immediate from the definition
\begin{lemma}\label{1} Let $\mathcal{C}^1_P=\{u\in \mathcal{C}^1(\R):u \hbox{ is $2\pi$ periodic}\} $. The set $\Psi=\{w \in \mathcal{C}^1_P: w'(x)\neq 0 \;\;\mbox{whenever}\;\; w(x)= 0 \}$ is an open dense subset of $\mathcal{C}^1_P$, $\ell(w)$ is finite if $w \in \Psi$ and $\ell$ is locally constant in $\Psi$.
\end{lemma}

The following result is due to Angenent (see \cite{Angenent}).
\begin{lemma}\label{angenent}
Let $q(t,x)$ and $r(t,x)$ be locally bounded functions in $(\tau, T)\times (-\pi,\pi)$ with $q_x$, $q_t$ locally bounded, and $w(t,x)$ be a classical solution of
\begin{equation}\label{eqlinear}
\begin{split}
&w_t= w_{xx} +q(t,x)w_x + r(t,x)w,\;\;x \in (-\pi,\pi), \;\;t\in(\tau, T)\\
&w(-\pi,t)=w(\pi,t),\ w_x(-\pi,t)=w_x(\pi,t), \quad t\geq 0.
\end{split}
\end{equation}
Suppose that $w$ is not identically zero. Then,
\begin{itemize}
  \item[(i)] $\ell(w(t, \cdot))$ is finite for each $t \in (\tau, T)$ and is monotone non-increasing in $t$.
  \item[(ii)] For each $t^* \in (\tau, T)$, $w(t, \cdot)$ belongs to $\Psi$ for each $t \in [t^* , T)$ except possibly for a finite number of points $t_1$, $\cdots$, $t_k$.
  \item[(iii)] If $w(t^* , \cdot) \notin \Psi$ for some $t^* \in (\tau, T)$, then
  $$\ell(w(t,\cdot))> \ell(w(s ,\cdot))$$
  for any $t \in (\tau , t^*)$ and $s \in (t^* , T)$.
\end{itemize}
\end{lemma}

\end{description}

\subsection{The construction of a positive global bounded solution when $\lambda >1$}
\label{maximal solution}

In this section we will derive a characterization of the pullback attractor $\{\mathscr{A}_\beta^\mathcal{C}(t):t\in \R\}$ of $\{T_\beta^\mathcal{C}(t,s): (t,s)\in \mathcal{P}\}$ where $T_\beta^\mathcal{C}(t,s)=T_\beta(t,s){\big|}_\mathcal{C}$ and $\mathcal{C}$ represents the positive cone within $H^1_0(0,\pi)$. The following lemmas will be helpful to obtain such characterization.

\begin{lemma} If $\mathcal{C}=\{\phi \in H^1_0(0,\pi): \phi(x)\geq 0, \;x\in [0,\pi]\}$, then $T_\beta(t,s)\mathcal{C}\subset \mathcal{C}$, for all $t\geq s$ and, if $0\neq u_0\in \mathcal{C}$, then
\begin{equation*}
\begin{split}
&u(t,s,u_0)(x)>0, \ x\in (0,\pi), \ t>s\quad \hbox{and}\\
&u_x(t,s,u_0)(0)>0, u_x(t,s,u_0)(\pi)<0, \ t>s.
\end{split}
\end{equation*}
\end{lemma}

\noindent{\bf Proof:} The result follows immediately from the Lap Number (Lemma \ref{angenent}) property of solutions of \eqref{naci} (after an odd $2\pi$-periodic extension of $u$) and the fact that $T_\beta(t,s)\mathcal{C}\subset \mathcal{C}$ for all $(t,s)\in \mathcal{P}$.\cqd

\begin{lemma} If $\lambda>1$ and $0\neq u_0\in \mathcal{C}$, then $S_\beta(t)u_0 \longrightarrow \phi_{1,\beta}^+$, as $t \to +\infty$.
\end{lemma}

\noindent{\bf Proof:} Recall that, for any $u_0\in H^1_0(0,\pi)$, $S_\beta(t)u_0\stackrel{t\to\infty}{\longrightarrow} \phi$ for some $\phi\in \mathcal{E}$, $S_\beta(t)\mathcal{C}\subset\mathcal{C}$ for all $t\geq 0$ and $\phi_{1,\beta}^+$ is the only element of $\mathcal{E}$ in $\mathcal{C}$. The result now follows trivially.\cqd

\medskip

Using these lemmas we can construct a global bounded solution in the positive cone which will play the role of the equilibria $\phi_{1,\beta}^+$ of \eqref{aci}; that is,
\begin{theorem}\label{existenceofxi}
If $\lambda>1$, $\{T_\beta(t,s): t \geq s\}$ has a global solution $\xi_1^+ : \R \to \mathcal{C}$ such that:

\begin{itemize}
\item [(i)] $\phi_{1,\beta_2}^+ \leq \xi_1^+ (t)\leq \phi_{1,\beta_1}^+$, for all $t \in \R$,
\item[(ii)] If  $0\neq u_0\in \mathcal{C}$ and $u_0\geq \xi_1^+ (s)$ for all $s\in \R$,  then $$\xi_1^+ (t)=\lim_{s\to -\infty} T_\beta(t,s) u_0.$$
\item[(iii)] Any bounded global solution $\psi:\R\to H^1_0(\Omega)$ of $\{T_\beta(t,s): t \geq s\}$ must satisfy $\psi(t)\leq \xi_1^+ (t)$, for all $t \in \R$.
\item[(iv)] If $\beta$ is constant $\xi_1^+ (t)=\phi_{1,\beta}^+$ for all $t\in \R$.
\end{itemize}
\end{theorem}

\noindent{\bf Proof:}  (i) Note that
\begin{equation*}
\begin{split}
\phi_{1,\beta_2}^+ = S_{\beta_2}(t-s)\phi_{1,\beta_2}^+ &\leq T_\beta(t,s) \phi_{1,\beta_2}^+
\\
&\leq T_\beta(t,s) \phi_{1,\beta_1}^+ \leq S_{\beta_1}(t-s)\phi_{1,\beta_1}^+=\phi_{1,\beta_1}^+,
\end{split}
\end{equation*}
since $\phi_{1,\beta_2}^+ \leq \phi_{1,\beta_1}^+$. Hence, for $s_1\leq s_2\leq t$,
\begin{equation*}
\begin{split}
\phi_{1,\beta_2}^+  &\leq T_\beta(t,s_1) \phi_{1,\beta_1}^+ = T_\beta(t,s_2) T_\beta(s_2,s_1)\phi_{1,\beta_1}^+ \\
&\leq T_\beta(t,s_2) S_{\beta_1}(s_2-s_1) \phi_{1,\beta_1}^+=  T_\beta(t,s_2) \phi_{1,\beta_1}^+\leq \phi_{1,\beta_1}^+.
\end{split}
\end{equation*}

Let $\xi_1^+ (t) =\lim_{s\to -\infty} T_\beta(t,s)\phi_{1,\beta_1}^+$.

\medskip

It is easy to see that $\xi_1^+ :\R\to H^1_0(0,\pi)$ is a global bounded solution of $\{T_\beta(t,s):t\geq s\}$ and that $\phi_{1,\beta_2}^+\leq \xi_1^+ (t) \leq \phi_{1,\beta_1}^+$ for all $t\in \R$.

(ii) If $0\neq u_0\in \mathcal{C}$, $u_0\geq \xi_1^+ (s)$ for all $s\in \R$, then $T_\beta(r,s)u_0\geq T_\beta(r,s)\xi_1^+ (s)=\xi_1^+ (r)$ for all $s\leq r$, so
\begin{equation*}
\begin{split}
\xi_1^+ (t)= T_\beta(t,r)\xi_1^+ (r)&
\leq T_\beta(t,r)\liminf_{s\to -\infty}T_\beta(r,s)u_0\\
&\leq T_\beta(t,r)\limsup_{s\to -\infty}T_\beta(r,s)u_0\\
&\leq
T_\beta(t,r)\lim_{s\to -\infty}S_{\beta_1}(r-s)u_0=T_\beta(t,r)\phi_{1,\beta_1}^+\stackrel{r\to-\infty}{\longrightarrow}\xi_1^+ (t),
\end{split}
\end{equation*}
and consequently $\xi_1^+ (t)=\lim_{s\to -\infty}T_\beta(t,s)u_0$.

(iii)  Let $\phi\in \mathcal{C}$ be such that $\phi\geq \psi(s)$ for all $s\in \R$ and $s_n\stackrel{n\to \infty}{\longrightarrow} -\infty$. Then,
\begin{equation*}
\begin{split}
\psi(t) &=T_\beta(t,s_n)\psi(s_n)\leq T_\beta(t,s_n)\phi\\
& \leq T_\beta(t,r)S_{\beta_1}(r-s_n)\phi \stackrel{n\to\infty}{\longrightarrow} T_\beta(t,r)\phi_{1,\beta_1}^+
\stackrel{r\to-\infty}{\longrightarrow}\xi_1^+ (t).
\end{split}
\end{equation*}
$\cqd$

To conclude this section we recall the characterization result for the pullback attractor of \eqref{naci} in the positive cone (see \cite{RB-VL}).

\begin{definition}\label{nondegeneracy}
A function $u:(-\infty,\tau]\to \mathcal{C}$ is said to be non-degenerate as $t\to -\infty$ if there exists $t_0\leq \tau$ and $\phi \in \mathcal{C}$ with $\phi(x)>0$ for all $x\in (0,\pi)$, $\phi'(0)>0$ and $\phi'(\pi)<0$ such that $u(t)\geq \phi$ for all $t\leq t_0$. Similarly we define non-degeneracy as $t\to +\infty$.
\end{definition}

\begin{theorem}[Rodriguez-Bernal \& Vidal-Lopez \cite{RB-VL}]
\label{Rodrigues-Bernal}
If $\lambda>1$, the global solution $\xi_1^+ :\R\to \mathcal{C}$ of $\{T_\beta(t,s): t \geq s \}$ given by Theorem \ref{existenceofxi} is the unique solution non-degenerate as $t\to - \infty$. If $0\neq u_0\in \mathcal{C}$, then $\|T_\beta(t,s)u_0-\xi_1^+ (t)\|_{H^1_0(0,\pi)} \stackrel{t\to +\infty}{\longrightarrow} 0$.
\end{theorem}

\section{Characterization of the $\omega$-limit and $\alpha$-limit sets}\label{S2}

In this section we study the $\omega$-limit sets of solutions and the $\alpha$-limit sets of global bounded solutions of \eqref{naci}.

We will study the $\alpha$-limit set adapting the ideas of  Chen and Matano in \cite{chen} for the $\omega$-limit set. The computations for the $\omega$-limit are analogous. To that end we first consider the same equation but with periodic boundary conditions
\begin{equation}\label{PP}
\left\{
\begin{array}{l}
  u_t= u_{xx} +\lambda u - \beta(t)u^3,\quad   x \in (-\pi,\pi),\;\;t>0, \\
  u(-\pi,t)=u(\pi,t),\ u_x(-\pi,t)=u_x(\pi,t), \quad t\geq 0,\\
  u(0)=u_0 \in H^1_P(-\pi,\pi).
\end{array}
\right.
\end{equation}
We know that the problem \eqref{PP} is globally well posed in $H^1_P(-\pi,\pi)$ and that the associated evolution process has a pullback attractor.

\begin{definition}
For each $2\pi$-periodic function $w$ and each $a \in \R$ we define the function $(\rho_{a}w)(x)= w(2a-x)$, $x\in \R$. The operator $w\mapsto \rho_{a}w$ is called \textbf{reflection}.
\end{definition}

\begin{remark}
The solutions $u(t,x)$ of (\ref{PP}) satisfy an equation of the type \eqref{eqlinear} with $q(t,x)=0$ and $r(t,x)= \lambda -\beta(t)u(t,x)^2$. Furthermore if $u$ also denotes the $2\pi$-periodic extension of $u$ to $\R$, the functions $\rho_{a}u$ restricted to $[-\pi,\pi]$ are also solutions of \eqref{PP} for each $a \in \R$.
\end{remark}

\begin{definition}
Let $v \in \mathcal{C}^1_P $. We say that $v$ is a \textbf{symmetrically oscillating function} if there exists a $x_0\in \R$ and $m \in \mathbb{N}$ such that
 \begin{eqnarray*}
  v(x) &=&v(2x_0 - x),\;\;x \in \R,  \\
  v'(x)&>& 0, \;\;x \in (x_0, x_0+ \pi/m), \\
  v(x)&=& v(x+ (2\pi)/m),\;\;x \in \R.
\end{eqnarray*}
We will denote the set of symmetrically oscillating functions by $\mathfrak{F}_m(x_0)$.
\end{definition}

\begin{figure}[h!]
\newrgbcolor{xdxdff}{0.49019607843137253 0.49019607843137253 1.}
\psset{xunit=1.0cm,yunit=1.0cm,algebraic=true,dimen=middle,dotstyle=o,dotsize=5pt 0,linewidth=1.6pt,arrowsize=3pt 2,arrowinset=0.25}
\begin{pspicture*}(-4,-1)(10,5)
\psplot[linewidth=0.8pt,plotpoints=200,linecolor=red]{-1.5}{6.2}{SIN(2.0*x)+3.0}
\psplot[linewidth=0.8pt]{-1.5}{6.4}{(-0.-0.*x)/3.141592653589793}
\psline[linewidth=0.8pt,linestyle=dashed,dash=5pt 5pt](3.9269908169872414,0.)(3.9269908169872414,4.)
\psline[linewidth=0.8pt,linestyle=dashed,dash=5pt 5pt](0.7853981633974483,0.)(0.7853981633974483,4.)
\rput[tl](2.2,-0.24){$x_0$}
\psline[linewidth=0.8pt,linestyle=dashed,dash=5pt 5pt](2.356194490192345,0.)(2.356194490192345,5.88)
\rput[tl](3.36,-0.14){$x_0\!+\!\frac{\pi }{m}$}
\rput[tl](0.1,-0.14){$x_0\!-\!\frac{\pi }{m}$}
\rput[tl](5.9,0.38){$\mathbb{R}$}
\begin{scriptsize}
\psdots[dotsize=5pt 0,dotstyle=*,linecolor=xdxdff](2.356194490192345,2.)
\psdots[dotsize=5pt 0,dotstyle=*,linecolor=xdxdff](3.9269908169872414,4.)
\psdots[dotsize=5pt 0,dotstyle=*,linecolor=xdxdff](0.7853981633974483,4.)
\end{scriptsize}
\end{pspicture*}
\caption{A common spatial symmetry property to the family of functions $\mathfrak{F}_m(x_0)$}
\end{figure}
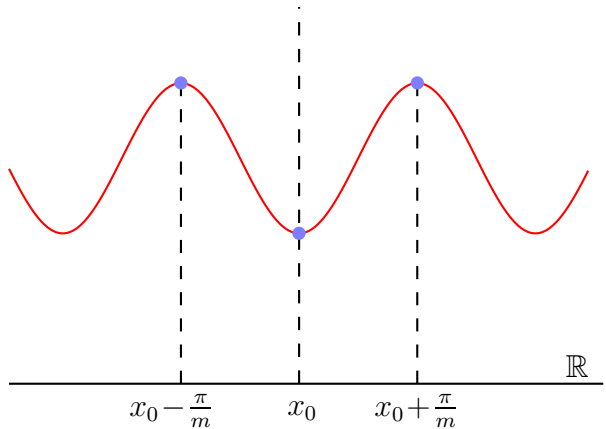



\begin{definition}
  Let $v \in \mathcal{C}^1([0,\pi])$ be such that $v(0)=v(\pi)=0$.  We say  $v$ is a \textbf{symmetrically oscillating function under the Dirichlet boundary conditions} if the odd $2\pi$-periodic extension of  $v$ belongs to either $\mathfrak{F}_m^+:=\mathfrak{F}_m(-\frac{\pi}{2m})$ or $\mathfrak{F}_m^-:=\mathfrak{F}_m(\frac{\pi}{2m})$ for some $m \in \mathbb{N}$. The set of all functions satisfying one of the above conditions is denoted by $\mathfrak{F}_m^\pm$.
\end{definition}

The main results in this section are the following:
\begin{theorem}\label{wlimite}
Let $u \in \mathcal{C}([0,+\infty), H^1(0,\pi))$ be the solution of problem
 \begin{equation}\label{PD2}
\left\{
\begin{array}{ll}
  u_t= u_{xx} +\lambda u - \beta(t)u^3, & 0<x<\pi,\;\;t>0, \\
  u(0,t)=u(\pi , t)=0, &\\
  u(0)=u_0. &
\end{array}
\right.
\end{equation}
Then there exists $m \in \mathbb{N}$ such that $\omega(u_0) \subset \mathfrak{F}_m^\pm \cup \{0\}$.
\end{theorem}

\begin{theorem}\label{alphalimite}
Let $\xi: \mathbb{R}\to H^1(0,\pi)$ be the global bounded  solution of problem
 \begin{equation}\label{PD2}
\left\{
\begin{array}{ll}
  u_t= u_{xx} +\lambda u - \beta(t)u^3, & 0<x<\pi,\;\;t\in \mathbb{R} \\
  u(0,t)=u(\pi , t)=0, \\
  u(0,x)=u_0(x) &
\end{array}
\right.
\end{equation}
Then there exists $m \in \mathbb{N}$ such that $\alpha_{\xi}(u_0) \subset \mathfrak{F}_m^\pm \cup \{0\}$.
\end{theorem}
\noindent
Before we prove the second theorem, we need some notations and preliminary lemmas.

  Let $\{t_n\}_{n\in \mathbb{N}}$ be a sequence in $\mathbb{R}$. For each $n \in \mathbb{N}$, let $\beta_n : \mathbb{R}\to \mathbb{R}$ be the function defined by $\beta_n(t)= \beta(t+ t_n)$. Under the assumptions of the function $\beta$, we have that the family $\{\beta_n\}_{n \in \mathbb{N}}$ is uniformly bounded and uniformly equicontinuous. Consequently, it has a subsequence (that we denote the same) and a globally Lipschitz and bounded function $\gamma : \mathbb{R} \to (0, +\infty)$ such that $\beta_n(t) \to \gamma(t)$ as $n \to +\infty$ uniformly in compact subsets of $\mathbb{R}$.

Now we consider the sequence of nonlinear problems
\begin{equation}\label{p0s1}
\left\{
\begin{array}{ll}
  u_t= u_{xx} +\lambda u - \beta(t)u^3,\quad   x \in (-\pi,\pi),\;\;t>0, \\
  u(-\pi,t)=u(\pi,t),\ u_x(-\pi,t)=u_x(\pi,t), \quad t\geq 0,\\
  u(0)=u_0 \in H^1_P(-\pi,\pi),
  \end{array}
\right.
\end{equation}
\begin{equation}\label{pns1}
\left\{
\begin{array}{ll}
 u_t= u_{xx} +\lambda u - \beta_n(t)u^3,\quad   x \in (-\pi,\pi),\;\;t>0, \\
  u(-\pi,t)=u(\pi,t),\ u_x(-\pi,t)=u_x(\pi,t), \quad t\geq 0,\\
  u(0)=u_0 \in H^1_P(-\pi,\pi),
\end{array}
\right.
\end{equation}
\begin{equation}\label{pinfinitos1}
\left\{
\begin{array}{ll}
  u_t= u_{xx} +\lambda u - \gamma(t)u^3,\quad   x \in (-\pi,\pi),\;\;t>0, \\
  u(-\pi,t)=u(\pi,t),\ u_x(-\pi,t)=u_x(\pi,t), \quad t\geq 0,\\
  u(0)=u_0 \in H^1_P(-\pi,\pi).
\end{array}
\right.
\end{equation}

Denote by $T_{\beta}(t,s)$, $T_{\beta_n}(t,s)$ and $T_{\infty}(t,s)$ the processes associated with (\ref{p0s1})-(\ref{pinfinitos1}) in $H^1_P(-\pi,\pi)=\{u\in H^1(-\pi,\pi): u(-\pi)=u(\pi)\}$.
We have that $T_{\beta}(t+t_n,s+t_n) = T_{\beta_n}(t,s)$  for all $t\geq s$. In fact, by the variation of constants formula
\begin{eqnarray*}
T_{\beta}(t+t_n,s+t_n)u_0 &=& e^{-A(t-s)}u_0 + \\
&&\int_{s+t_n}^{t+t_n}e^{-A(t+t_n - \theta)}\Big\{\lambda T_{\beta}(\theta,s+t_n)u_0 - \beta(\theta)\Big[T_{\beta}(\theta,s+t_n)u_0\Big]^3\Big\}d\theta\\
&=& e^{-A(t-s)}u_0 + \\
&&\int_{s}^{t}e^{-A(t- \bar{\theta})}\Big\{\lambda T_{\beta}(\bar{\theta} + t_n,s+t_n)u_0 - \beta_n(\bar{\theta})\Big[T_{\beta}(\bar{\theta} + t_n,s+t_n)u_0\Big]^3\Big\}d\bar{\theta}.
\end{eqnarray*}
Since $T_{\beta}(t+t_n,s+t_n)u_0$ and $T_{\beta_n}(t,s)u_0$ are solutions of the same integral equation we have the result.

 Now, let $\xi: \mathbb{R}\to H^1_P(-\pi,\pi)$ be a global bounded solution of (\ref{p0s1}). Define the sequence of functions $\xi_n: \mathbb{R} \to H^1_P(-\pi,\pi)$ by $\xi_n(t)= \xi(t+ t_n)$. We prove that this sequence has a subsequence which converges to a global bounded solution of problem (\ref{pinfinitos1}). To simplify the notation, we define the functions
\begin{equation}
\label{functions}
F(t,u)=\lambda u - \beta(t) u^3,\;\;\;\; F_n(t,u)= \lambda u -\beta_n(t)u^3,\;\;\;\;F_{\infty}(t,u)=\lambda u - \gamma(t)u^3,
\end{equation}
 for $(t,u)\in \mathbb{R}\times H^1_P(-\pi,\pi)$. With the hypotheses on the function $\beta$ we can easily prove that the functions $F$, $F_n$ and $F_{\infty}$ take bounded subsets of $H^1_P(-\pi,\pi)$ in uniformly bounded subsets (in $t$) of $H^1_P(-\pi,\pi)$. With this, we can prove that  $\xi$ is also bounded in $X^1: = \{u \in H^2(-\pi,\pi)\cap H^1_P(-\pi,\pi): u'(-\pi)=u'(\pi)\}$. In fact, as the linear semigroup decays exponentially we have $\xi$ is a solution of the integral equation

\begin{equation}
  \xi(t)= \int_{-\infty}^t e^{-A(t-s)}F(s,\xi(s))ds.
\end{equation}
Then,
\begin{eqnarray*}
  \|\xi(t)\|_1& \leq & \int_{-\infty}^t\|Ae^{-A(t-s)}F(s,\xi(s))\|_{L^2}ds \\
  &\leq & \int_{-\infty}^t\|A^{1/2}e^{-A(t-s)}\|\;\|F(s,\xi(s))\|_{1/2}ds \\
  &\leq& \int_{-\infty}^t C_{1/2}(t-s)^{1/2}e^{-\delta(t-s)}\;K ds\\
  &= &C_{1/2}K\delta^{-1/2}\Gamma(1/2) < +\infty,
\end{eqnarray*}
 where we use that $\sup_{s\in\mathbb{R}}\|F(s,\xi(s))\|_{1/2} \leq K$, for some constant $K$ and recall that $X^{1/2}= H^1_P(-\pi,\pi)$.
Therefore, $$\sup_{t\in\mathbb{R}}\{\|\xi(t)\|_{1/2},\;\|\xi(t)\|_1,\; \|\xi_t(t)\|_{L^2}\} \;<\;\infty.$$
Thus, by the Arzel\`a-Ascoli Theorem, we have that the sequence  $\xi_n$ in $\mathcal{C}(\mathbb{R}, H^1_P(-\pi,\pi))$
 has a subsequence which converges uniformly in compact subsets of $\R$ to a continuous function $\zeta: \mathbb{R} \to H^1_P(-\pi,\pi)$.
Now, as $\xi_n(t)= \xi(t + t_n)$ we have

\begin{eqnarray*}
\xi_n(t)&=& \int_{-\infty}^{t+t_n}e^{-A(t+t_n-s)}[\lambda \xi(s)-\beta(s)\xi(s)^3]ds \\
  &=& \int_{-\infty}^{t}e^{-A(t-s)}[\lambda \xi(s+t_n)-\beta(s+t_n)\xi(s+t_n)^3]ds \\
  &=& \int_{-\infty}^{t}e^{-A(t-s)}[\lambda \xi_n(s)-\beta_n(s)\xi_n(s)^3]ds.
\end{eqnarray*}
From this, it is not difficult to see that
\begin{equation*}
  \zeta(t)= \int_{-\infty}^t e^{-A(t-s)} [\lambda \zeta(s)- \gamma(s)\zeta(s)^3]ds
\end{equation*}
and, in particular, $\zeta$ is a global bounded solution of the problem (\ref{pinfinitos1}).

The following lemma plays a fundamental role and is adapted from Lemma 3.7 \cite{chen}.
\begin{lemma}\label{primeiro}
Let $\xi : \mathbb{R} \to H^1_P(-\pi,\pi)$ be a global bounded solution of
\begin{equation}\label{PP2}
\left\{
\begin{array}{ll}
  u_t= u_{xx} +\lambda u - \beta(t)u^3,\quad   x \in (-\pi,\pi),\;\;t>0, \\
  u(-\pi,t)=u(\pi,t),\ u_x(-\pi,t)=u_x(\pi,t), \quad t\geq 0,\\
  u(0)=u_0 \in H^1_P(-\pi,\pi),
\end{array}
\right.
\end{equation}
and let  $\varphi$ be an element of the set $\alpha_{\xi}(u_0)$. Then, for each $a \in [-\pi, \pi]$, we have
\begin{itemize}
  \item[(i)] either $\rho_{a}\varphi = \varphi$ or $\rho_{a}\varphi - \varphi \in \Psi$,
  \item[(ii)] $\rho_{a}\varphi = \varphi$ if and if  $\varphi'(a)=0$.
\end{itemize}
\end{lemma}
 \noindent{\bf Proof:} We need only to prove the assertion (i), because (ii) follows immediately from (i). To prove (i), take a sequence $t_n \to + \infty$ such that $\xi(-t_n) \to \varphi$ in $C^1_P$. Define $\beta_n(t)= \beta (t - t_n)$, $t \in \mathbb{R}$. The family $\{\beta_n\}_{n \in \mathbb{N}}$ is uniformly bounded and uniformly  equicontinuous. Consequently, it has a subsequence (which we denote the same) and a globally Lipschitz and bounded function  $\gamma : \mathbb{R} \to [\beta_1, \beta_2]$ such that $\beta_n(t) \to \gamma(t)$ as $n \to +\infty$ uniformly in compact subsets of $\mathbb{R}$. Similarly, we define $\xi_n(t)= \xi(t-t_n)$, for all $t \in \mathbb{R}$ and $n \in \mathbb{N}$.

 Since $\sup_{t\in \mathbb{R}}\{\|\xi(t)\|_{H^1},\;\|\xi(t)\|_1,\; \|\xi_t(t)\|_{L^2}\} < \infty$ there is a subsequence,  which again we denote by $\{\xi_n\}$, converging to a function, denoted by $p$, which satisfies
\begin{equation}\label{pns12}
\left\{
\begin{array}{l}
  p_t= p_{xx} +\lambda p - \gamma(t)p^3,  x \in (-\pi,\pi),\;\;t \in \mathbb{R}, \\
  p(-\pi,t)=p(\pi,t), \ p_x(-\pi,t)=p_x(\pi,t),\;\;t \in \mathbb{R},\\
  p(0,x)= \varphi(x), \;\;x \in (-\pi,\pi).
\end{array}
\right.
\end{equation}
Now, we define the function $w = \rho_{a}p-p$, which is a solution of linear problem
\begin{equation*}
\left\{
\begin{array}{ll}
  w_t =w_{xx} + r_{\infty}(t,x)w, \;\; x \in (-\pi,\pi),\;\;t \in \mathbb{R},& \\
  w(-\pi,t)=w(\pi,t), \ w_x(-\pi,t)=w_x(\pi,t),\;\;t \in \mathbb{R},&\\
  w(0,x)= (\rho_{a}\varphi - \varphi)(x), \;\;x \in (-\pi,\pi),&
\end{array}
\right.
\end{equation*}
where  $r_{\infty}(t,x)= \lambda - \gamma (t)\Big[\frac{(\rho_{a}p)^3 - p^3}{\rho_{a}p-p}\Big]$.
 Suppose that $\rho_{a}\varphi \neq \varphi$. From Lemma \ref{angenent}, there is $\delta >0$ such that $w(\delta,\cdot)$ and $w(-\delta,\cdot)$ have only simple zeroes, that is, belong to $\Psi$. But
$$(\rho_{a}\xi_n - \xi_n)(\pm \delta,\cdot )\;=\; (\rho_{a}\xi - \xi)(\pm \delta - t_n,\cdot) \to w(\pm \delta, \cdot)$$
in $C^1_P$. So, there exists a positive integer $N_1$ such that
\begin{equation}\label{2}
\ell\Big((\rho_{a}\xi -\xi)(\delta - t_n,\cdot)\Big)\;=\; \ell (w(\delta,\cdot)),
\end{equation}
 for all $n \geq N_1$. But $\rho_a\xi-\xi$ satisfies a parabolic equation of the form (\ref{eqlinear}) with $q\equiv 0$ and $r$ locally bounded, so we conclude from  Lemma \ref{angenent}(i)  and (\ref{2}) that
 \begin{equation}\label{3}
  \ell\Big((\rho_{a}\xi -\xi)(t,\cdot)\Big)\;=\; \ell (w(\delta,\cdot)),
\end{equation}
for all $t\leq \delta - t_{N_1}$. Similarly, choosing a positive integer $N_2$ sufficiently large, we have
\begin{equation}\label{4}
  \ell\Big((\rho_{a}\xi -\xi)(t,\cdot)\Big)\;=\; \ell (w(-\delta,\cdot)).
\end{equation}
 for all $t \leq -\delta - t_{N_2}$. It follows from (\ref{3}) and (\ref{4}) that $\ell (w(\delta,\cdot))=\ell (w(-\delta,\cdot))$. From Lemma \ref{angenent}(iii), this implies that $w(0,\cdot)= \rho_a\varphi-\varphi \in \Psi$. The proof of item (i) is complete.$\cqd$

 \begin{lemma}\label{segundo}
Let $\xi$ be as in Lemma \ref{primeiro}. Then the set
$$\alpha_{\xi}(u_0) \subset \bigcup_{
\begin{array}{c}x\in \R\\
  m \in \mathbb{N}
   \end{array}
} \mathfrak{F}_m(x)\cup \mathfrak{F},$$
where $\mathfrak{F}$ is the set of all constant functions.
\end{lemma}
 \noindent{\bf Proof:} Let $\varphi$ be a nonconstant element of set $\alpha_{\xi}(u_0)$. Choose $x^*$ such that $\varphi'(x^*)\neq 0$. Without loss of generality , suppose that $\varphi'(x^*) > 0$ and let $I =(x_0, x_1)$ be the  maximal interval containing $x^*$ such that $\varphi'(x) >0$ for all $x \in I$. Since $\varphi$ is $\mathcal{C}^1$, we have that $\varphi'(x_0)= \varphi'(x_1)=0$. It follows from Lemma \ref{primeiro} that $\rho_{x_0}\varphi = \varphi = \rho_{x_1}\varphi$ implying that there is $m \in \mathbb{N}$ such that $x_1-x_0 = \pi/m$. Thus $\varphi \in \mathfrak{F}_m(x_0)$.$\cqd$

\bigskip

\noindent{\bf Proof of Theorem \ref{alphalimite}:}
Let $\varphi$ and $\psi$ be nonzero functions belonging to the set $\alpha_{\xi}(u_0)$, $\{t_n\}_{n\in \mathbb{N}}$ and $\{s_n\}_{n \in \mathbb{N}}$ sequences tending to $+\infty$ such that $\xi(-t_n) \to \varphi$ and $\xi(-s_n)\to \psi$ in $C^1([0,\pi])$. Let $\tilde{\xi}$, $\tilde{\varphi}$ and $\tilde{\psi}$ be the odd extensions  $2\pi$-periodic of functions  $\xi$, $\varphi$ and $\psi$. Then $\tilde{\xi}(t,x)$ is a global bounded solution of a  nonautonomous nonlinear parabolic equation, with periodic boundary conditions on the interval $[-\pi ,\pi]$. Furthermore, $\tilde{\xi}(-t_n)\to \tilde{\varphi}$ and  $\tilde{\xi}(-s_n) \to \tilde{\psi}$ in  $C^1([-\pi,\pi])$. From Lemma \ref{segundo}, there are $x_0 \in [-\pi, \pi]$ and  $m \in \mathbb{N}$ such that $\tilde{\varphi} \in \mathfrak{F}_m(x_0)$. Note that $x_0 \neq 0$ since  $\tilde{\varphi}(-x) \neq \tilde{\varphi}(x)$, this implies that $\tilde{\varphi}'(0) \neq 0$. It follows from this, Lemma \ref{primeiro} and \ref{segundo} that  $\tilde{\varphi}$ has only simple zeroes. Let us prove that $\tilde{\psi}$ also belongs to the set $\mathfrak{F}_m(x_0)$. Since the convergence of $\tilde{\xi}(-t_n)$ to $\tilde{\varphi}$ is in $C^1_P$ and $\tilde{\varphi}$ has only simple zeroes, we ensure the existence  of $n_0 \in \mathbb{N}$ such that $\tilde{\xi}(-t_n) \in \Psi$  and
$$\ell(\tilde{\xi}(-t_n))\;=\; \ell (\tilde{\varphi}),$$
for all $n \geq n_0$. But $\tilde{\xi}(t)$ is a solution of a parabolic equation of the form (\ref{eqlinear}), then $\ell (\tilde{\xi}(t))$ is a non-increasing function of $t$. So,
\begin{equation}\label{conclusao1}
\ell (\tilde{\xi}(-t))\;=\; \ell (\tilde{\varphi})
\end{equation}
for all $t \geq t_{n_0}$. Similarly, we can conclude that there is a positive integer $n_1$ such that
\begin{equation}\label{conclusao1}
\ell (\tilde{\xi}(-t))\;=\; \ell (\tilde{\psi})
\end{equation}
 for all $t \geq s_{n_1}$. Therefore, for all $t \geq \max \{t_{n_0}, s_{n_1}\}$, we have that $\ell (\tilde{\xi}(-t))\;=\; \ell (\tilde{\varphi})\;=\; \ell(\tilde{\psi})$ and $\tilde{\xi}(-t) \in \Psi$.

 Since the functions of the set $\mathfrak{F}_m(x_0)$ have $2\pi/m$ as the fundamental period, we have that  $\tilde{\varphi}'(x) >0$ for all $x_0<x<x_0+\frac{1}{2}\Big(\frac{2\pi}{m}\Big)$. Without loss of generality, we suppose  that $\tilde{\varphi}'(0) >0$. Since $\tilde{\varphi}$ is odd and $\frac{2\pi}{m}$-periodic, we conclude that
$$
\tilde{\varphi}'(x) > 0,\;\;\mbox{for all}\;\;x \in  \Big(-\frac{\pi}{2m},  \frac{\pi}{2m}\Big),
$$
consequently, $\tilde{\varphi} \in \mathfrak{F}_m\Big(-\frac{\pi}{2m}\Big)$, in other words, $\varphi \in \mathfrak{F}_m^+$. Since $\tilde\varphi'(0) >0$ we conclude that $\tilde{\xi}_x(-t_n,0) >0$ for all $n \geq N_1$, for some $N_1$. So, we can see that $\tilde{\psi}'(0)> 0$, as otherwise  $\tilde{\psi}'(0)< 0$ and therefore $\tilde{\xi}_x(-s_n,0) <0$, for all $n \geq N_2$, for some $N_2$. This would imply the existence of $t_n^*$ between $t_n$ and $s_n$ such that $\tilde{\xi}_x(-t_n^*,0)=0$, which would contradict the fact that $\tilde{\xi}(-t) \in \Psi$ for all $t \geq \max \{t_{n_0}, s_{n_1}\}$. Thus, $\tilde{\psi}'(0)> 0$ and since $\ell (\tilde{\varphi})= \ell(\tilde{\psi})$ we conclude that $\tilde{\psi} \in \mathfrak{F}_m(x_0)$, where $x_0 = -\frac{\pi}{2m}$. Note that assuming $\tilde{\varphi}'(0) <0$, $\tilde{\varphi}$ and $\tilde{\psi}$ would belong to  $\mathfrak{F}_m\Big(\frac{\pi}{2m}\Big)$, i.e., $\varphi$ and $\psi$ would belong to the set $\mathfrak{F}_m^-$. $\cqd$

\section{The structure of the pullback attractor for $ 0 < \lambda < 4$}
\label{S3}

 In \cite{Alexandre} the authors proved that if $0 < \lambda \leq 1$, then all solutions of problem \eqref{naci} tend to the trivial solution. One of our results in this paper, was to determine the structure of the pullback attractor when the parameter $1 < \lambda < 4$. We know from  \cite{Alexandre} that the pullback attractor is the set of  all global bounded solutions  of problem \eqref{naci}. In Section \ref{maximal solution} we proved the existence of two global bounded solutions, $\pm \xi_1^+ (t)$, which are non-degenerate as $t \to \pm \infty$. Let us find the others global bounded solutions when $1 < \lambda < 4 $. For this end, we resume to the analysis of the $\alpha$-limit sets and $\omega$-limit sets. In principle, we consider the solutions which  $w$-limit set and $\alpha$-limit set are different from the unitary set $\{0\}$. According to Theorems \ref{wlimite} and \ref{alphalimite} we have the following possible cases:

{\bf $1{\bf.}\!\!^{\rm st}$ Case: $\omega(u_0) \subset \mathfrak{F}_1^+ \cup \{0\}$}\\
  The proof of this case is valid for any $\lambda > 1$. Indeed, assume there is a nonzero function $\varphi \in \omega(u_0)$ and, that is, a sequence $t_n \to +\infty$ such that $T_\beta(t_n,0)u_0 = u(t_n, \cdot) \to \varphi$ in $C^1([0,\pi])$ as $n \to +\infty$. Furthermore, $\varphi(x) >0$ for all $0<x<\pi$ and $\rho_{\frac{\pi}{2}}\varphi = \varphi$. Thus, for some positive integer  $n_0$ sufficiently large  we have
  $$  \Big[T_\beta(t_{n_0},0)u_0\Big](x) = u(t_{n_0}, x) \;>\;0,$$
  for all $0<x<\pi$,  so $0\neq u(t_{n_0})\in \mathcal{C}$ and from Theorem \ref{Rodrigues-Bernal} we can conclude $$\|T_\beta(t,t_{n_0})u(t_{n_0})-\xi_1^+ (t)\|_{H^1_0(0,\pi)} \stackrel{t\to +\infty}{\longrightarrow} 0.$$
 Then, $0 \not \in \omega(u_0)$ and  $\omega(u_0)\subset \omega(\xi^+_1(0))\subset [\phi_{\beta_2}^+, \phi_{\beta_1}^+]$.

Analogously, if $\omega(u_0) \subset \mathfrak{F}_1^- \cup \{0\}$ we conclude that $0 \not \in \omega(u_0)$ and $\omega\left(  u_{0}\right)\subset \omega(-\xi^+_1(0))\subset [\phi_{\beta_1}^-, \phi_{\beta_2}^-]$. $\cqd$\\

{\bf $2{\bf.}\!\!^{\rm nd}$ Case: $\alpha_{\xi}(u_0) \subset \mathfrak{F}_1^+ \cup \{0\}$}\\
  As in the first case, this proof is valid for any $\lambda > 1$. Assume that there is a nonzero  function  $\psi \in \alpha_{\xi}(u_0) \subset \mathfrak{F}_1^+$, that is, a sequence $t_n \to +\infty$  such that $\xi(-t_n) \to \psi$ in $C^1([0,\pi])$. Again, $\psi(x) >0$ for all $0<x<\pi$ and then, by Lemma \ref{angenent},
  $\xi(t,x) \geq 0$, for all $t \in \mathbb{R}$ and $0 < x < \pi$. Let us prove that $\xi$ is non-degenerate as $t \to - \infty$, i.e., there exists $t^* \in \mathbb{R}$ and a nonzero function   $\vartheta \geq 0$  such that $\xi(-t) \geq \vartheta$ for all $t> t^*$.  In fact, consider $\tilde{\beta}_2$ sufficiently large, $\beta(t)\leq \beta_2 \leq \tilde{\beta}_2$, so that the positive equilibrium for the autonomous problem with $\tilde{\beta}_2$ instead of $\beta(t)$ satisfies:
\begin{equation*}
    \phi_{1,\tilde{\beta}_2}^+ \leq \frac{1}{2}\psi .
  \end{equation*}
  So, there exists $n_0$ such that, for all $n \geq n_0$, we have
    \begin{equation*}
    \phi_{1,\tilde{\beta}_2}^+ \leq \frac{1}{2}\psi \leq \xi(-t_n) .
  \end{equation*}
   Let us see that $\phi_{1,\tilde{\beta}_2}^+ \leq \xi(-t)$ for all $t \geq t_{n_0}$. Fixed $t \geq t_{n_0}$ consider $m \geq n_0$ such that $t_{m} \leq t < t_{m+1}$. From the comparison results, we have
  \begin{equation*}
    \phi_{1,\tilde{\beta}_2}^+ = S_{\tilde{\beta}_2}(t_{m+1}-t)\phi_{1,\tilde{\beta}_2}^+\leq T_{\beta}(-t,-t_{m+1})\phi_{1,\tilde{\beta}_2}^+\leq T_{\beta}(-t,-t_{m+1})\xi(-t_{m+1}) = \xi(-t),
  \end{equation*}
  as we wanted. Hence, from Theorem \ref{Rodrigues-Bernal}, we conclude that  $\xi(t)=\xi_1^+ (t)$, since $\xi_1^+ $ is the unique positive solution which is non-degenerate as $t \to -\infty$. Thus, $0 \notin \alpha_{\xi}(u_0)=\alpha_{\xi_1^+ }(u_0)$ and $\alpha_{\xi}(u_0) \subset [\phi_{\beta_2}^+, \phi_{\beta_1}^+]$.

Analogously, if $\alpha_{\xi}(u_0) \subset \mathfrak{F}_1^- \cup \{0\}$ we conclude that $\xi(t)= -\xi_1^+ (t)$ and consequently $0 \notin \alpha_{\xi}(u_0)=\alpha_{-\xi_1^+ }(u_0)$ and $\alpha_{\xi}(u_0) \subset [\phi_{\beta_1}^-, \phi_{\beta_2}^-]$. $\cqd$\\

\begin{remark} The second case, which we have just seen, consider global bounded positive solutions which do not tend to trivial solution as $t \to - \infty$. As a corollary of this proof, we have the uniqueness of solution for $t$ negative for such initial data. In fact, if $\xi_1(t)$ and $\xi_2(t)$ are global solutions which pass through $u_0$ such that $\alpha_{\xi_1}(u_0)$ and $\alpha_{\xi_2}(u_0)$ are contained in the set $\mathfrak{F}_1^+ \cup \{0\}$, it follows that both solutions are non-degenerate as $t\to -\infty$ and, from the results of \cite{Alexandre}, we have $\xi_1(t)=\xi_2(t)=\xi_1^+ (t)$.
\end{remark}

{\bf $3{\bf.}\!\!^{\rm rd}$ Case: $\omega(u_0) \subset \mathfrak{F}_j^{\pm}\cup \{0\}$, $j=2,3,\cdots$.}\\
    Let us show that these cases do not occur for nonzero functions when $\lambda \in (1,4)$. To simplify the proof we will make it for $\mathfrak{F}_2^+$, the proof is analogous for the other cases. Suppose by contradiction that there exists a non-zero function $\varphi \in \omega (u_0)$ and consider a sequence $t_n \to + \infty $ such that $ T_{\beta}(t_n,0)u_0  =  u (t_n, \cdot) \to \varphi $ in $ C^1([0, \pi])$ as $ n \to +\infty $. As we mentioned, $\varphi \in \mathfrak{F}_2^+$ implies that $\varphi(x)>0$ for $0<x<\pi/2$ and $\varphi(\pi -x)= - \varphi(x)$, for all $x \in [0,\pi]$.

     Since $T_{\beta}(t,0)u_0= u(t,\cdot)$ is bounded, we may define $u_n(t,\cdot)= u(t+t_n,\cdot)$, for all $t\geq -t_n$, and $\beta_n(t)= \beta(t+t_n)$, $t \in \mathbb{R}$. Then there are subsequences of $\{u_n\}_n$ and $\{\beta_n\}_n$ which converge, respectively, to $p$ e $\gamma$, uniformly in compact subsets  of $\mathbb{R}$, being $p$ a global bounded solution of the initial value problem
   \begin{equation}\label{pns122}
\left\{
\begin{array}{ll}
  p_t= p_{xx} +\lambda p - \gamma(t)p^3, & 0<x<\pi,\;\;t \in \mathbb{R}, \\
  p(0,t) = p(\pi,t)=0, \\
  p(0,x)= \varphi(x). &
\end{array}
\right.
\end{equation}
 Since $\varphi(\pi -x)= - \varphi(x)$, for all $x \in [0,\pi]$, and from the uniqueness solution we have $p(t,\pi-x)= -p(t,x)$, for all $x \in [0,\pi]$ and $t \geq 0$. Consequently, $p(t,\pi/2)=0$ for all $t \geq 0$. Let us see that $p(t,\pi/2)=0$ also for $t<0$. In fact, fix $t>0$ and consider $t_{n_0}$ such that $t \leq t_{n_0}$. So, for all $n \geq n_0$, we have $u_n(-t,\cdot)$ well defined and $u_n(-t,\cdot)= u(t_n-t, \cdot) \to p(-t, \cdot)$. Then, $p(-t, \cdot) \in \omega(u_0)$ and, consequently, in this case, $p(-t,\cdot) \in \mathfrak{F}_2^+$ implying $p(-t,\pi/2)=0$. In this way , the equation in (\ref{pns122}) can be considered, together with the Dirichlet boundary conditions, in the interval  $[0,\pi/2]$. However, the operator $-\frac{\partial^2}{\partial x^2}$ with the Dirichlet boundary conditions in the interval  $[0,\pi/2]$ has the eigenvalues  $\lambda_n= 4n^2$, $n \in \mathbb{N}$. Since the parameter $\lambda \in (1,4)$, we are in the case where, according \cite{Alexandre}, any solution tend to zero. Thus, the unique global bounded solution of equation in (\ref{pns122}) with Dirichlet boundary conditions in the interval $[0,\pi/2]$ is the null solution, contradicting the fact that $p(t,\cdot)$ is a global bounded nonzero solution. Thus, for $\lambda \in (1,4)$, the set $\omega(u_0)$ cannot be contained in the set $\mathfrak{F}_2^+$. $\cqd$ \\

{\bf $4{\bf.}\!\!^{\rm th}$ Case: $\alpha_{\xi}(u_0) \subset \mathfrak{F}_j^{\pm} \cup \{0\}$, $j=2,3,\cdots $.}\\
For the same reasons given in the previous case we have that such cases cannot occur for nonzero functions for $\lambda \in (1,4)$ . $\cqd$\\

Now, for $N^2 < \lambda < (N+1)^2$, $N \geq 1$, we prove that there is no non-constant solution connecting the zero equilibrium to itself (homoclinic orbit). Suppose by contradiction that $\xi(t) \rightarrow 0$ as $t \rightarrow \pm \infty$. Consider the linearized operator about the trivial solution $L = \partial_x^2 +\lambda I_d$. We know that the eigenvalues are defined by $\lambda_n=\lambda -n^2$, $n=1,2,\cdots$, and the eigenfunctions $\psi_n(x)=\sin(nx)$ have $n+1$ simple zeros in $[0,\pi]$. In this way, zero is not an eigenvalue and $\lambda_1 >\cdots >\lambda_N >0 >\lambda_{N+1}> \cdots >\lambda_n \rightarrow -\infty$ implying that the trivial equilibrium is hyperbolic and, from \cite{Henry}, its local unstable manifold is an N-dimensional Lipschitz manifold and the rate of approach from $\xi(t)$ to zero is exponential. Furthemore  $\xi(t)$ is a nontrivial solution of the linear parabolic equation
$$v_t = v_{xx} + \lambda v -\beta(t)\xi(t,x)^2\ v, \;\; 0<x<\pi ,\;\;v(0,t)=v(\pi,t)=0.$$
Define the operator $L(t): D(L(t))\subset L^2(0,\pi) \rightarrow L^2(0,\pi)$, $L(t)v = v_{xx}+ \lambda v -\beta(t)\xi(t,x)^2\ v$, where $D(L(t))= H^2(0,\pi)\cap H^1_0(0,\pi)$. If $\lambda_n(t)$, $n =1,2,\cdots$, are the eigenvalues of $L(t)$, note that $\lambda_n(t)\rightarrow \lambda_n$ as $t \rightarrow \pm \infty$. Now, we use some results on asymptotic behavior of solutions of linear parabolic equations from \cite{Henry2} (see Theorem 3, 4 and 5). From \cite[Theorem 4]{Henry2} we conclude that there is an integer $j\geq 1$ and a constant $C_1\neq 0$ such that
\begin{equation*}
\xi(t,x)= {\rm exp}\left(\int_0^t\lambda_j(s)ds\right)[C_1\psi_j(x)+ o(1)]
\ \hbox{ as }\ t \rightarrow \infty,
\end{equation*}
with convergence in the sense of $C^1([0,\pi])$. Since $\xi(t)\rightarrow 0$, $\psi_j$ has to be related to a negative eigenvalue, that is, $j\geq N+1$. From $C^1$ convergence, we conclude that $\xi(t)$ vanishes $j+1$ times in $[0,\pi]$ for $t$ large enough.
 On the other hand, from \cite[Theorems 3 and 5]{Henry2}, we also conclude that there is an integer $k\geq 1$ and a constant $C_2\neq 0$ such that
\begin{equation*}
\xi(t,x)= exp\left(-\int_t^0\lambda_k(s)ds\right)[C_2\psi_k(x)+ o(1)]\ \hbox{\  as \ } t \rightarrow -\infty,
\end{equation*}
with convergence in the sense of $C^1([0,\pi])$. Again, since $\xi(t)\rightarrow 0$, $\psi_k$ has to be related a positive eigenvalue, that is, $k\leq N$. From $C^1$ convergence, we conclude that $\xi(t)$ vanishes $k+1$ times in $[0,\pi]$ for $-t$ large enough.
Finally, applying Lemma \ref{angenent} we see that $j+1\leq k+1$, consequently, $N+2\leq j+1\leq k+1\leq N+1$, a contradiction, proving that there is no homoclinic orbit at zero.

Consequently, in the case $1< \lambda < 4$, that is, $N=1$, the unstable manifold of the trivial solution is $1$-dimensional, then a global bounded solution that goes out from zero tends to one of the solutions $\pm \xi^+_1$ as $t \rightarrow  \infty$, according to the analysis of the $\omega$-limit sets seen in this section. Without loss of generality, if the solution tends to $\xi^+_1(t)$ as $t \rightarrow  \infty$, it must be positive in $(0,\pi)$ for $t$ sufficiently large. But from Lemma \ref{angenent}, we conclude the solution is positive for all $t \in \mathbb{R}$
(see Proposition \ref{lambda<4} for more detailed explanation).

We summarize the results of this section in the next theorem.

 \begin{theorem}\label{lambdaentre1e4}
 If $1< \lambda < 4$, the unstable manifold of the zero is one-dimensional and the global bounded solutions, $\pm \zeta(t)$,  which leave from the zero, tend to the non-degenerate solutions $\pm \xi_1^+ (t)$. In particular, the pullback attractor is given by $\mathcal{A}(t)=W^u(0)(t)\cup\{\xi_1^+(t)\}\cup\{-\xi^+_1(t)\}, \;\; t \in \mathbb{R}$. In addition, for any $u_0\in H^1_0(0,\pi)$, $\omega(u_0)\subset \omega(\xi_1^+ (0))\cup\omega(-\xi_1^+ (0))\cup \{0\}$.
 \end{theorem}

This theorem, besides giving a thorough characterization of the pulback attractor $\{\mathcal{A}(t), t\in \R\}$, also shows that the pullback and the skew product attractor associated to \eqref{naci} have gradient structure.

\section{ The structure of the attractor for $N^2 < \lambda < (N+1)^2$, $N \geq 2$}
\label{valero}
We have seen in the two previous sections that the  non-degenerate global bounded solutions $\pm \xi_1^+ (t):= \xi_1^{\pm}(t)$ of  \eqref{naci} play the same role, when $1<\lambda < 4$, that the equilibria $\phi_{1,\beta}^{\pm}$ of the autonomous problem. Now, assuming that $N^2 < \lambda < (N+1)^2$ and knowing that there are $2N$ non-zero equilibria of the autonomous equation, $\phi_{j,\beta_i}^{\pm}$, $j= 1, \cdots , N$ and $i=1,2$, we can construct the following global bounded solutions of \eqref{naci}:
\begin{equation}\label{equinaoauto}
  \xi_j^{\pm}(t)=\lim_{s\to - \infty} T_{\beta}(t,s)\phi_{j,\beta_1}^{\pm},
\end{equation}
$j= 2, \cdots , N$.

\begin{figure}[h!]
\psset{xunit=1.0cm,yunit=1.0cm,algebraic=true,dimen=middle,dotstyle=o,dotsize=10pt 0,linewidth=1.6pt,arrowsize=3pt 2,arrowinset=0.25}
\begin{pspicture*}(-0.7,-2)(8,3.4)
\psplot[linewidth=0.6pt,plotpoints=200,linecolor=blue]{0}{6.28}{0.5*SIN(x)}
\psplot[linewidth=0.6pt,linestyle=dashed,dash=2pt 2pt,plotpoints=200,linecolor=red]{0}{6.28}{0.8*SIN(x)}
\psplot[linewidth=0.6pt,plotpoints=200,linecolor=blue]{0}{6.28}{1.5*SIN(x)}
\psline[linewidth=0.8pt](0,-1.5)(0,2)
\psplot[linewidth=0.8pt]{-0.5}{7}{(-0.-0.*x)/7.854}
\begin{scriptsize}
\rput[tl](6.3,-0.07){$\pi$}
\rput[tl](2.95,-0.07){$\frac{\pi}{2}$}
\rput[tl](1.3,1.9){$\phi^+_{2,\beta_{1}}$}
\rput[tl](1.3,0.5){$\phi_{2,\beta_2}^+$}
\rput[tl](1.3,1.2){$\xi_2^+ (t)$}
\rput[tl](-0.2,-0.07){$0$}
\end{scriptsize}

\end{pspicture*}
\caption{The second ``non-autonomous equilibrium"}\label{xi_2^+}
\end{figure}
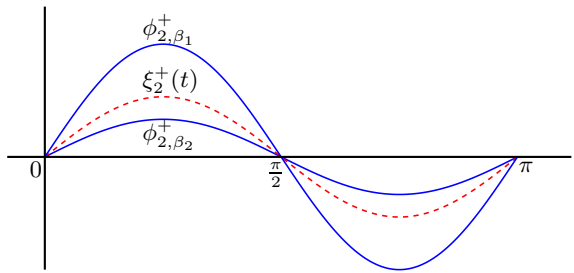

Just as we did for $\xi_1^+ (t)$, we can show that each solution $\xi_j^\pm(t)$ is non-degenerate as $t \to \pm \infty$, belongs to the sets $\mathfrak{F}_j^\pm$ and
$$Y_j^\pm= \Big\{v \in H_0^1(0,\pi): \min(\phi_{j,\beta_1}^\pm(x), \phi_{j,\beta_2}^\pm(x))\leq v(x)\leq \max(\phi_{j,\beta_1}^\pm(x), \phi_{j,\beta_2}^\pm(x))\Big\}.$$
Furthermore, it is not difficult to prove that if $u_0 \in \mathfrak{F}_j^\pm$ then $\omega (u_0) \subset\omega(\xi_j^\pm(0))$. We will call such solutions by ``non-autonomous equilibria" .

\begin{lemma}
\label{ExistB0} Assume that $N^2<\lambda<(N+1)^2$, $N\geq 2$ and that $\omega ( u_0)  \subset \mathfrak{F}_{j}^{\pm} \cup \{0\}$, $2\leq j \leq N$. If $0\neq \varphi \in \omega (u_0)$ then, there exists $\widetilde{\varphi}\in \omega(  u_{0} )  \cap Y_{j}^{\pm}$.
\end{lemma}

 \noindent{\bf Proof:} We prove the case when $j$ is even and $\varphi \in \mathfrak{F}_j^+$. There exists a sequence $t_{n}\rightarrow+\infty$ such that
$T_{\beta}(  t_{n},0)  u_{0}=u(  t_{n})  \rightarrow \varphi$ in $C^{1}([0,\pi])$. Moreover, $\varphi$ satisfies

\begin{align}
 \varphi \left(\frac{\pi}{j}-x\right)  &=\varphi(x) \text{ if } x \in \left[0,\frac{\pi}{j}\right],\; \varphi \left(\frac{\pi}{j}\right)=0, \nonumber\\
 \varphi \left(\frac{2\pi}{j}-x\right)  &= - \varphi(x) \text{ if } x \in \left[0,\frac{2\pi}{j}\right], \; \varphi\left(\frac{2\pi}{j}\right)=0, \nonumber\\
 \varphi \left(\frac{3\pi}{j}-x\right)  &= \varphi(x) \text{ if } x \in \left[0,\frac{3\pi}{j}\right], \; \varphi \left(\frac{3\pi}{j}\right)=0, \nonumber\\
 &\ \vdots \label{PropFi}\\
  \varphi \left(\pi-x\right)  &= - \varphi(x) \text{ if } x \in [0,\pi],\nonumber\\
\varphi(x)  & >0, \text{ if } x \in  \mathrm{I_1} \doteq \bigcup_{k=0}^{\frac{j-2}{2}}\left(\frac{2k\pi}{j},\frac{(2k+1)\pi}{j}\right),\nonumber\\
\varphi(x)  & <0, \text{ if } x  \in  \mathrm{I_2} \doteq \bigcup_{k=0}^{\frac{j-2}{2}}\left(\frac{(2k+1)\pi}{j},\frac{(2k+2)\pi}{j}\right).\nonumber
\end{align}
The figure below sketches the function $\varphi$ when j is even.  If $j$ is odd,  $\varphi(\pi-x)=\varphi(x)$ for all $ x\in [0,\pi]$ and, consequently, $\varphi(x)>0$, $\frac{(j-1)\pi}{j}<x<\pi$.

\begin{figure}[h!]
\newrgbcolor{xdxdff}{0.5 0.5 1}
\psset{xunit=3.0cm,yunit=3.0cm,algebraic=true,dimen=middle,dotstyle=o,dotsize=5pt 0,linewidth=0.2pt,arrowsize=3pt 2,arrowinset=0.25}
\begin{pspicture*}(-1.2,-0.5)(6,1)
\psplot[linewidth=0.6pt,plotpoints=200,linecolor=red]{0}{1.4}{SIN(6.0*x)/3.0}
\psplot[linewidth=0.6pt,plotpoints=200,linecolor=red]{1.7}{3.15}{SIN(6.0*x)/3.0}
\psline[linewidth=0.6pt](0,-0.5)(0,0.8)
\psline[linewidth=0.6pt](-0.5,0)(3.4,0)
\begin{scriptsize}{\color{blue}
\rput[tl](1.4,0.274){$.$}
\rput[tl](1.42,0.246){$.$}
\rput[tl](1.44,0.216){$.$}
\rput[tl](1.67,-0.198){$.$}
\rput[tl](1.65,-0.17){$.$}
\rput[tl](1.63,-0.14){$.$}}
\rput[tl](0.42,-0.03){$\frac{\pi}{j}$}
\rput[tl](1.05,-0.03){$\frac{2\pi}{j}$}
\rput[tl](0.3,0.42){$\varphi$}
\rput[tl](3.16,-0.03){$\pi$}
\rput[tl](2.30,-0.03){$\frac{(j-1)\pi}{j}$}
\rput[tl](-0.1,-0.03){$0$}
\psdots[dotsize=3pt 0,dotstyle=*,linecolor=xdxdff](0.5235987755982988,0.)
\psdots[dotsize=3pt 0,dotstyle=*,linecolor=xdxdff](1.0471975511965976,0.)
\psdots[dotsize=3pt 0,dotstyle=*,linecolor=xdxdff](3.141592653589793,0.)
\psdots[dotsize=3pt 0,dotstyle=*,linecolor=xdxdff](2.617993877991494,0.)
\psdots[dotsize=3pt 0,dotstyle=*,linecolor=xdxdff](0.,0.)
\end{scriptsize}
\end{pspicture*}
 \caption{Case $j$ even}\label{F_j^+}
\end{figure}
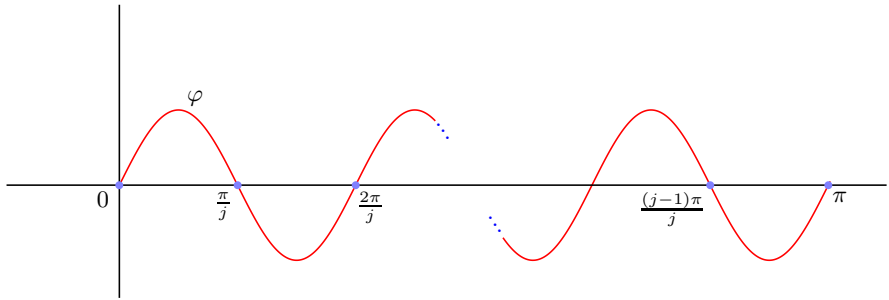


Let $v_{s}\left(  t\right)  =T_{\beta}\left(  t,s\right)  \varphi$, $t \geq s$. It is clear
that
\begin{equation}\label{PropVi}
v_{s}(t,x)  =\left\{
\begin{array}{ll}
v_{s}^{1}(t,x)  \text{ if }0\leq x\leq\frac{\pi}{j},\\
v_{s}^{2}(t,x)  \text{ if }\frac{\pi}{j}\leq x\leq \frac{2\pi}{j},\\
\qquad\vdots\\
v_{s}^{j}(t,x)  \text{ if }\frac{(j-1)\pi}{j}\leq x\leq \pi,
\end{array}
\right.
\end{equation}
where $v_{s}^{1}(t)$, $v_{s}^{2}(t)$, $\ldots$, $v_{s}^{j}(t)$ are the solutions of (\ref{naci})
with initial data $\varphi$ but restricted to the intervals $[0,\frac{\pi}{j}]$, $[\frac{\pi}{j},\frac{2\pi}{j}]$, $\ldots$,
$[\frac{(j-1)\pi}{j},\pi]$, respectively. Denote by $\varphi^{1}$, $\varphi^{2}$, $\ldots$, $\varphi ^j$ the restriction
of $\varphi$ to the intervals $[0,\frac{\pi}{j}]$, $[\frac{\pi}{j},\frac{2\pi}{j}]$, $\ldots$,
$[\frac{(j-1)\pi}{j},\pi]$, respectively. Also, $S_{\beta_{i}}^{1}$, $S_{\beta_{i}}^{2}$, $\ldots$, $S_{\beta_{i}}^{j}$ will be the
corresponding semigroups for $\beta$ constant in the same intervals. Then by
comparison we have
\begin{equation}\label{scomp_vi}
\begin{array}{lllll}
0\leq   & S_{\beta_{2}}^{1}(t-s)\varphi^{1}     &\leq v_{s}^{1}(t)   &\leq S_{\beta_{1}}^{1}(t-s)\varphi^{1},    & \\
        & S_{\beta_{1}}^{2}(t-s)\varphi^{2}     &\leq v_{s}^{2}(t)   &\leq S_{\beta_{2}}^{2}(t-s)\varphi^{2}     & \leq 0,\\
        & \qquad \vdots                         &    \qquad \vdots   &     \qquad \vdots                         & \\
0\leq   & S_{\beta_{1}}^{j-1}(t-s)\varphi^{j-1} &\leq v_{s}^{j-1}(t) &\leq S_{\beta_{2}}^{j-1}(t-s)\varphi^{j-1},& \\
        & S_{\beta_{1}}^{j}(t-s)\varphi^{j}     &\leq v_{s}^{j}(t)   &\leq S_{\beta_{2}}^{j}(t-s)\varphi^{j}     &\leq 0,
\end{array}
\end{equation}
assuming $j$ is even. If $j$ is odd, the last expression would be
$$0\;\leq \;S_{\beta_{2}}^{j}(t-s)\varphi^{j}\; \leq\; v_{s}^{j}(t) \; \leq \;S_{\beta_{1}}^{j}(t-s)  \varphi^{j},$$
since $\varphi^j(x) >0$, for all $\frac{(j-1)}{j}<x<\pi$.

Recall that the eigenvalues of the operator $-\partial_x^2$ with Dirichlet conditions on the interval $\left[0,\frac{\pi}{j}\right]$ are $\lambda_{n,\frac{\pi}{j}}= n^2j^2$, $n=1,2,\cdots$. Then
$\lambda > N^2 \geq j^2 = \lambda_{1,\frac{\pi}{j}}$. Hence, for any $\varepsilon>0$, there exists $T(\varepsilon) >0$ such that
\begin{align*}
\left\| S_{\beta_{1}}^{1}(t-s)  \varphi^{1}-\phi_{1,\beta_{1}
,\frac{\pi}{j}}^{+}\right\| _{L^\infty\left(  0,\frac{\pi}{j}\right)  } &
\leq \varepsilon,\\
\left\| S_{\beta_{2}}^{1}(t-s)  \varphi^{1}-\phi_{1,\beta_{2}
,\frac{\pi}{j}}^{+}\right\| _{L^\infty\left(  0,\frac{\pi}{j}\right)  } &
\leq \varepsilon,
\end{align*}
if $t-s\geq T(\varepsilon)$, where $\phi_{1,\beta_{1},\frac{\pi}{j}}^{+}$  and $\phi_{1,\beta_{2}
,\frac{\pi}{j}}^{+}$ represent the equilibria of the semigroups $S_{\beta_{1}}^{1}$, $S_{\beta_{2}}^{1}$ in $H_0^1\left(0, \frac{\pi}{j}\right)$, respectively. Hence

\begin{equation}
\label{Ineq1}
\phi_{1,\beta_{2},\frac{\pi}{j}}^{+}(x) - \varepsilon \leq v_{s}^{1}\left(
t,x\right)  \leq\phi_{1,\beta_{1},\frac{\pi}{j}}^{+}(x) + \varepsilon
\end{equation}
if $t-s\geq T(\varepsilon)$. In the same way,

\begin{equation}\label{Ineq2}
 -\phi_{1,\beta_{1},\frac{\pi}{j}}^{+}\left(\frac{2\pi}{j}-x\right)-\varepsilon   \leq v_{s}^{2}(t,x)\leq
 -\phi_{1,\beta_{2},\frac{\pi}{j}}^{+} \left(\frac{2\pi}{j}-x\right)
+\varepsilon,
\end{equation}
$$\vdots$$
\begin{equation}\label{Ineq22}
 -\phi_{1,\beta_{1},\frac{\pi}{j}}^{+}(\pi-x)-\varepsilon  \leq v_{s}^{j}(t,x)\leq
 -\phi_{1,\beta_{2},\frac{\pi}{j}}^{+}(\pi-x)
+\varepsilon
\end{equation}
if $t-s\geq T(\varepsilon)$.

Further, we can show that there is a constant $R>0$, which does not depend on $t$ and $s$ (see Proposition 12.8 \cite{CLR12}), such that
\begin{align*}
\left\Vert u\left(  t\right)  \right\Vert _{H_{0}^{1}(0,\pi)}  &  \leq R,\ \forall
t\geq 0,\\
\left\Vert v_{s}\left(  t\right)  \right\Vert _{H_{0}^{1}(0,\pi)}  &  \leq R,\text{
}\forall t\geq s.
\end{align*}
Using the variation of constants formula, Gronwal's inequality and the embedding of $H^1$ into $L^\infty$, the difference satisfies
\begin{equation}\label{difference}
\|u(t)-v_{s}(t)\|_{L^\infty}\leq e^{\delta (t-s)}\|u(s)-\varphi\|_{H_{0}^{1}}
\end{equation}
for some $\delta = \delta(R)>0$. Since $u(t_{n})\rightarrow \varphi$ in $H_0^1(0,\pi)$, for any $\varepsilon>0$, $T(\varepsilon)$ (where
$T(\varepsilon)$ is taken from (\ref{Ineq1})) there exists $t_{n_{\varepsilon}}$ such that
$$
\|u(t_{n_{\varepsilon}})-\varphi\|_{H^1_0} \leq \varepsilon e^{-\delta T(\varepsilon)}.
$$
Hence,
\begin{equation}
\|u(t_{n_{\varepsilon}}+T(\varepsilon)) -v_{t_{n_{\varepsilon}}}(t_{n_{\varepsilon}}+T(\varepsilon))\|_{L^\infty} \leq \varepsilon. \label{Ineq3}
\end{equation}
It follows from (\ref{Ineq1})-(\ref{Ineq3}) that
\begin{align*}
\phi_{j,\beta_{2}}^{+}(x)  -2\varepsilon &  \leq u(t_{n_{\varepsilon}}+T(\varepsilon),x) \leq \phi
_{j,\beta_{1}}^{+}(x)  + 2\varepsilon,
\end{align*}
for all $x \in \mathrm{I}_1$ and
\begin{align*}
\phi_{j,\beta_{1}}^{+}(x)  -2\varepsilon &  \leq u(t_{n_{\varepsilon}}+T(\varepsilon),x)  \leq \phi_{j,\beta_{2}}^{+}(x)  +2\varepsilon,
\end{align*}
for all $x \in \mathrm{I}_2$.

We choose $\varepsilon_{m}\rightarrow 0$. Then, passing to a subsequence
$$u(t_{n_{\varepsilon_m}}+T(\varepsilon_m))\rightarrow\widetilde{\varphi}\in \mathfrak{F}_{j}^{+} \text{ in }C^{1}([0,\pi]),
$$
and
\begin{align*}
\phi_{j,\beta_{2}}^{+}(x)   &  \leq \widetilde{\varphi}(x) \leq \phi_{j,\beta_{1}}^{+}(x),
\end{align*}
for all $x \in \mathrm{I}_1$, and
\begin{align*}
\phi_{j,\beta_{1}}^{+}(x)   &  \leq \widetilde{\varphi}(x) \leq \phi_{j,\beta_{2}}^{+}(x),
\end{align*}
for all $x \in \mathrm{I}_2$. Therefore, $\widetilde{\varphi} \in \omega(u_{0})\cap Y_j^+$.$\cqd$

\bigskip

Let us define the hull of $\beta(\cdot)$ by
\begin{equation}\label{hull}
\mathcal{H}(\beta)  =cl_{C(\mathbb{R},\mathbb{R})} \{\beta(\cdot + s):\;s \in \mathbb{R}\},
\end{equation}
where the closure is taken with respect to the metric of the uniform convergence in compact subsets of $\R$. Also, let $\xi_{j,\gamma}^{\pm}(\cdot)$ be the ``non-autonomous equilibria" in $\mathfrak{F}_{j}^{\pm}$ for problem (\ref{naci}) with $\beta(\cdot)$ replaced by $\gamma(\cdot)$.

\begin{theorem}\label{charac-omega}
Let $N^2<\lambda<(N+1)^2$, $N \geq 2$, $\omega(u_{0}) \subset \mathfrak{F}_j^{\pm} \cup \{0\}$ and let $\varphi \in \omega(u_{0})$
be such that $\varphi \neq 0$. Then $\varphi \in Y_{j}^{\pm}$. Hence, $\omega(u_{0})  \subset Y_{j}^{\pm}$. Moreover,
\begin{equation}
\omega(u_{0}) \subset \{\xi_{j,\gamma}^{\pm}(t):\gamma \in \mathcal{S}^+(\beta),\; t \in \mathbb{R}\},
\label{FormulaOmega}
\end{equation}
where $\mathcal{S}^+(\beta)$ is given by Definition \ref{s-pm-beta}.
\end{theorem}

\noindent{\bf Proof:} 
We note that all sets involved in the asymptotics are compact in $C^1([0,\pi])$ and therefore we may use indistinctly $H^1(0,\pi)$ or $C^1([0,\pi])$ convergence.
Again, we prove the case for $\mathfrak{F}_{j}^{+}$. Suppose by contradiction that $\varphi \notin Y_{j}^{+}\cap\, C^1([0,\pi])$. Since
$Y_{j}^{+}$ is closed, there exist disjoint open neighbourhoods (in $C^1([0,\pi])$) $\mathcal{O}$, $\mathcal{O}_{\varphi}$ of $Y_{j}^{+}\cap\, C^1([0,\pi])$ and $\varphi$, respectively.

By Lemma \ref{ExistB0} there is $\widetilde{\varphi} \in \omega(u_{0})\cap Y_{j}^{+}\,\cap\, C^1[0,\pi]$. Then we can choose a sequence $t_{m}\rightarrow +\infty$ such that $u(t_{m})  \rightarrow \widetilde{\varphi}$, where $u(t)=T_{\beta}(t,0)u_{0}$.

We note that $u(t_{m})  \in \mathcal{O}$ and that for each $t_{m}$ there exists a first time $\sigma_{m}$ such that
\begin{align*}
u(t)   &  \in \mathcal{O} \text{ for }t_{m}\leq t<t_{m}+\sigma_{m},\\
u(t_{m} + \sigma_{m})   &  \in \overline{\mathcal{O}},\\
u(t)   &  \not \in \mathcal{O} \text{ for } t_{m} + \sigma_{m} \leq t \leq t_{m}+\sigma_{m}+T_{m},
\end{align*}
for some $T_m >0$.
The sequence $\sigma_{m}$ goes to $+\infty$. Indeed, we define $v_{m}(t)=u(t+t_{m})$, which is a solution of (\ref{naci})
with $\beta_{m}(t)=\beta(t+t_{m})$. Up to a subsequence, $v_{m}(t)$ converges uniformly in compact sets to a complete bounded solution $q(t)$ of problem (\ref{naci}) but replacing $\beta(t)$ by the limit of $\beta_{m}(t)$, denoted by $\gamma(t)$. It follows that $q(t) \in \omega(u_0)\subset \mathfrak{F}_{j}^{+}\cup\{0\}$, for all $t\in\mathbb{R}$, and that $q(0)=\widetilde{\varphi}$.
 Restricting the system to the subintervals $[0,\frac{\pi}{j}]$, $[\frac{\pi}{j},\frac{2\pi}{j}]$, $\cdots$, and using comparison, we obtain that $q(t)  \in Y_{j}^{+} \subset \mathcal{O}$, for all $t\geq 0$. If $\sigma_{m}$ is bounded, then we
can assume that $\sigma_{m} \rightarrow \sigma$, and $v_{m}(\sigma_{m})  \notin \mathcal{O}$ implies that $q(\sigma)\notin \mathcal{O}$,
which is a contradiction.

Further, we define the sequence $u_{m}(t)=u(t+t_{m}+\sigma_{m})$. It is clear that
\begin{align*}
u_{m}(t)   &  \in \mathcal{O}\text{ for }-\sigma_{m}\leq t<0,\\
u_{m}(0)   &  \in \overline{\mathcal{O}},\\
u_{m}(t)   &  \not \in \mathcal{O} \text{ for } 0\leq t\leq T_{m}.
\end{align*}
This function is a solution of problem (\ref{naci}) with $\overline{\beta}_{m}(t)=\beta(t+t_{m}+\sigma_{m})$.
Up to a subsequence, $u_{m}(t)$ converges uniformly in compact sets to a complete bounded solution $p(t)$ of problem (\ref{naci}) but replacing $\beta(t)$ by the limit of $\overline{\beta}_{m}(t)$, denoted by $\overline{\gamma}(t)$. It follows that $p(t)\in \mathfrak{F}_j^+\cap\{0\}$ for all $t$,  and that $p(t)  \in \overline{\mathcal{O}}$, for all $t\leq0$.


We note that there is a function $\phi\in C^{1}([0,\pi])$ such that
$(-1)^{k-1}p(t)\geq(-1)^{k-1}\phi$ in $[\frac{(k-1)\pi}{j},\frac{k\pi}{j}]$,
$1\leq k\leq j$, $t\leq0$, $(-1)^{k-1}\phi_{x}(\frac{(k-1)\pi}{j})>0,$ $1\leq
k\leq j+1,$ and $(-1)^{k-1}\phi(x)>0$, for $x\in\left(  \frac{(k-1)\pi}%
{j},\frac{k\pi}{j}\right)  $, $1\leq k\leq j$. This follows from the fact that
$p(\frac{k\pi}{j},t)=0$, $0\leq k\leq j$, $p(t)\in C^{1}([0,\pi])$ and that
$p(t)$ lies in a $C^{1}([0,\pi])$ small neighbourhood of $Y_{j}^{+}\cap
C^{1}([0,\pi])$ for all $t\leq0$.

%
%

Restricting the system to the subintervals $[0,\frac{\pi}{j}]$,
$[\frac{\pi}{j},\frac{2\pi}{j}]$, $\cdots$, we obtain that $p(t)$ is the unique
nondegenerate solution as $t \rightarrow -\infty$, so that $p(t)= \xi_{j, \bar{\gamma}}^+(t) \in Y_{j}^{+}\subset \mathcal{O}$ for all $t\in\mathbb{R}$.
Furthermore, $u_{m}(0) \rightarrow p(0)$ and $u_m(0) \notin \mathcal{O}$ imply that $p(0)\not \in Y_{j}^{+}$, a contradiction.
It follows that $\varphi \in Y_{j}^{+}$.

As $\omega(u_{0})$ is connected, we have that $0\notin \omega(u_{0})$, and then $\omega(u_{0}) \subset Y_{j}^{+}$.

Finally, let us prove (\ref{FormulaOmega}). Let $ \varphi \in \omega(u_{0})  \subset Y_{j}^{+}$. We can choose a sequence $t_{n}
\rightarrow +\infty$ such that $u(t_{n})  \rightarrow \varphi$, where $u(t)=T_{\beta}(t,0)u_{0}$. Arguing as before we
can prove that $u_{n}(t)=u(t+t_{n})$ converges to a complete bounded solution $p(t)$ of problem (\ref{naci}) but replacing $\beta(t)$
 by the limit of $\beta_n(t)= \beta(t+t_{n})$, denoted by $\gamma(t)  \in \mathcal{S}^+(\beta)$.
Since $p(t)  \in Y_{j}^{+}$ for all $t\in\mathbb{R}$, it is clear
that $p(t)=\xi_{j,\gamma}^{+}(t)$. Hence, $\varphi = p(0)  =\xi_{j,\gamma}^{+}(0)$.$\cqd$

\smallskip
The characterization of the $\alpha-$limit is obtained in a similar way with suitable changes.

\begin{lemma}
\label{ExistB0Alfa}Let $N^{2}<\lambda<(N+1)^{2}$, $N\geq2$, and let $\xi$ be a
bounded global solution with initial condition $u_{0}$ such that $\alpha_{\xi
}(u_{0})\subset\mathfrak{F}_{j}^{\pm}\cup\{0\}$, $2\leq j\leq N$. Suppose the
existence of $\varphi\in\alpha_{\xi}(u_{0})$ such that $\varphi\neq0$, then
there exists $\widetilde{\varphi}\in\alpha_{\xi}(u_{0})\cap Y_{j}^{\pm}$.
\end{lemma}

\noindent\textbf{Proof:} We prove the case $\mathfrak{F}_{j}^{+}$ with $j$
even. There exists a sequence $t_{n}\rightarrow+\infty$ such that $\xi
(-t_{n})\rightarrow\varphi$ in $C^{1}([0,\pi])$. Moreover, $\varphi$ satisfies \eqref{PropFi}.

If $j$ is odd, $\varphi(\pi-x)=\varphi(x)$ for all $x\in\lbrack0,\pi]$ and,
consequently, $\varphi(x)>0$, $\frac{(j-1)\pi}{j}<x<\pi$.

Let $v_{s}\left(  t\right)  =T_{\beta}\left(  t,s\right)  \varphi$, $t\geq s$.
It is clear that $v_s$ satisfies \eqref{PropVi} with $v_s^j$ and $\varphi^j$ are as before. Also with $S_{\beta_{1}}^{j}(t-s)\varphi^{j}$ as before, by comparison we have \eqref{scomp_vi} for $j$ even.

If $j$ is odd, the last expression would be%
\[
0\;\leq\;S_{\beta_{2}}^{j}(t-s)\varphi^{j}\;\leq\;v_{s}^{j}(t)\;\leq
\;S_{\beta_{1}}^{j}(t-s)\varphi^{j},
\]
since $\varphi^{j}(x)>0$, for all $\frac{(j-1)}{j}<x<\pi$.


As before, for any $\varepsilon>0$, there exists
$T(\varepsilon)>0$ such that.
\begin{align*}
\left\Vert S_{\beta_{1}}^{1}(t-s)\varphi^{1}-\phi_{1,\beta_{1},\frac{\pi}{j}%
}^{+}\right\Vert _{L^\infty\left(  0,\frac{\pi}{j}\right)  } &
\leq\varepsilon,\\
\left\Vert S_{\beta_{2}}^{1}(t-s)\varphi^{1}-\phi_{1,\beta_{2},\frac{\pi}{j}%
}^{+}\right\Vert _{L^\infty\left(  0,\frac{\pi}{j}\right)  } &
\leq\varepsilon,
\end{align*}
if $t-s\geq T(\varepsilon)$, where $\phi_{1,\beta_{1},\frac{\pi}{j}}^{+}$ and
$\phi_{1,\beta_{2},\frac{\pi}{j}}^{+}$ represent the equilibria of the
semigroups $S_{\beta_{1}}^{1}$, $S_{\beta_{2}}^{1}$, respectively. We can choose $T(\varepsilon)$ such
that $T(\varepsilon)\rightarrow+\infty$ as $\varepsilon\rightarrow0$.

Hence
\begin{equation}
\phi_{1,\beta_{2},\frac{\pi}{j}}^{+}(x)-\varepsilon\leq v_{s}^{1}\left(  x,t\right)  \leq\phi_{1,\beta
_{1},\frac{\pi}{j}}^{+}(x)+\varepsilon\label{Ineq1B}
\end{equation}
if $t-s\geq T(\varepsilon)$. In the same way,
\begin{equation}
-\phi_{1,\beta_{1},\frac{\pi}{j}}^{+}\left(  \frac{2\pi}{j}-x\right)
-\varepsilon\leq v_{s}^{2}(x,t)\leq
-\phi_{1,\beta_{2},\frac{\pi}{j}}^{+}\left(  \frac{2\pi}{j}-x\right)
+\varepsilon\label{Ineq2B}
\end{equation}

\[
\vdots
\]

\begin{equation}
-\phi_{1,\beta_{1},\frac{\pi}{j}}^{+}(\pi-x)-\varepsilon\leq v_{s}^{j}(x,t)\leq-\phi_{1,\beta_{2},\frac{\pi}{j}}%
^{+}(\pi-x)+\varepsilon\label{Ineq22B}%
\end{equation}
if $t-s\geq T(\varepsilon)$.

Further, we can show as before that there is a constant $R>0$, which does not
depend on $t$ and $s$ such that%
\begin{align*}
\left\Vert \xi\left(  t\right)  \right\Vert _{H_{0}^{1}(0,\pi)} &  \leq
R,\ \forall t\geq0,\\
\left\Vert v_{s}\left(  t\right)  \right\Vert _{H_{0}^{1}(0,\pi)} &  \leq
R,\text{ }\forall t\geq s.
\end{align*}
Proceeding as before we obtain \eqref{difference}, and for any $\varepsilon>0$ and
$T(\epsilon)$ taken from \eqref{Ineq1B} there exists $t_{n_{\varepsilon}}\geq 2T(\varepsilon)$ such that

\[
\Vert\xi(-t_{n_{\varepsilon}})-\varphi\Vert_{H^1_0}\leq\varepsilon
e^{-\delta T(\varepsilon)}.
\]
Hence,%
\begin{equation}
\Vert\xi(-t_{n_{\varepsilon}}+T(\varepsilon))-v_{-t_{n_{\varepsilon}}%
}(-t_{n_{\varepsilon}}+T(\varepsilon))\Vert_{L^{\infty}}\leq\varepsilon
.\label{Ineq3B}%
\end{equation}

It follows from (\ref{Ineq1B})-(\ref{Ineq3B}) that%
\[
\phi_{j,\beta_{2}}^{+}(x)-2\varepsilon  \leq\xi(x,-t_{n_{\varepsilon}}+T(\varepsilon))\leq\phi_{j,\beta_{1}%
}^{+}(x)+2\varepsilon ,
\]
for all $x\in\mathrm{I}_{1}$ and%
\[
\phi_{j,\beta_{1}}^{+}(x)-2\varepsilon  \leq\xi(x,-t_{n_{\varepsilon}}+T(\varepsilon))\leq\phi_{j,\beta_{2}%
}^{+}(x)+2\varepsilon  ,
\]
for all $x\in\mathrm{I}_{2}$.

We choose $\varepsilon_{m}\rightarrow0$. Then, as $-t_{n_{\varepsilon_{m}}%
}+T(\varepsilon_{m})\leq-2T(\varepsilon_{m})\rightarrow-\infty$, passing to a
subsequence%
\[
\xi(-t_{n_{\varepsilon_{m}}}+T(\varepsilon_{m}))\rightarrow\widetilde{\varphi
}\in\mathfrak{F}_{j}^{+}\text{ in }C^{1}([0,\pi]),
\]
and%
\[
\phi_{j,\beta_{2}}^{+}(x)\leq\widetilde{\varphi}(x)\leq\phi_{j,\beta_{1}}%
^{+}(x),\ \text{for all }x\in\mathrm{I}_{1},
\]%
\[
\phi_{j,\beta_{1}}^{+}(x)\leq\widetilde{\varphi}(x)\leq\phi_{j,\beta_{2}}%
^{+}(x),\ \text{for all }x\in\mathrm{I}_{2}.
\]
Therefore, $\widetilde{\varphi}\in\alpha_{\xi}(u_{0})\cap Y_{j}^{+}$.$\cqd$

\bigskip

Now we are ready to characterize the $\alpha-$limit of points in the phase space.


\begin{theorem}\label{charac-alpha}
Let $N^{2}<\lambda<(N+1)^{2}$, $N\geq2$, and let $\xi$ be a bounded global
solution with initial conditon $u_{0}$ such that $\alpha_{\xi}(u_{0}%
)\subset\mathfrak{F}_{j}^{\pm}\cup\{0\}$, $2\leq j\leq N$. Suppose the
existence of $\varphi\in\alpha_{\xi}(u_{0})$ be such that $\varphi\neq0$. Then
$\varphi\in Y_{j}^{\pm}$. Hence, $\alpha_{\xi}(u_{0})\subset Y_{j}^{\pm}$.
Moreover,%
\begin{equation}
\alpha_{\xi}(u_{0})\subset\{\xi_{j,\gamma}^{\pm}(t):\gamma\in\mathcal{S}^-%
(\beta),\;t\in\mathbb{R}\}.\label{FormulaAlpha}%
\end{equation}

\end{theorem}

\noindent\textbf{Proof:} 
We note that all sets involved in the asymptotics are compact in $C^1([0,\pi])$ and therefore we may use indistinctly $H^1(0,\pi)$ or $C^1([0,\pi])$ convergence.
Again, we prove the case for $\mathfrak{F}_{j}^{+}$.
Suppose by contradiction that $\varphi\notin Y_{j}^{+}$. Since $Y_{j}^{+}$ is
closed, there exist disjoint open neighbourhoods $\mathcal{O}$, $\mathcal{O}%
_{\varphi}$ of $Y_{j}^{+}$ and $\varphi$, respectively.

By Lemma \ref{ExistB0Alfa} there is $\widetilde{\varphi}\in\alpha_{\xi}%
(u_{0})\cap Y_{j}^{+}$. Then we can choose sequences $s_{n},\ t_{m}%
\rightarrow+\infty$ such that $\xi\left(  -s_{n}\right)  \rightarrow\varphi$,
$\xi(-t_{m})\rightarrow\widetilde{\varphi}$.

There exists $N>0$ such that $\xi(-t_{m})\in\mathcal{O}$, $\xi(-s_{n}%
)\in\mathcal{O}_{\varphi}$ for $n,m\geq N$. Hence, for each $n\geq N$ there
exist $t_{m_{n}}>s_{n}$ and $\sigma_{n}, T_{n}>0$ such that%
\begin{align*}
\xi(t) &  \in\mathcal{O}\text{ for }-t_{m_{n}}\leq t<-t_{m_{n}}+\sigma_{n},\\
\xi(-t_{m_{n}}+\sigma_{n}) &  \in\overline{\mathcal{O}},\\
\xi(t) &  \not \in \mathcal{O}\text{ for }-t_{m_{n}}+\sigma_{n}\leq
t\leq-t_{m_{n}}+\sigma_{n}+T_{n}.
\end{align*}

The sequence $\sigma_{n}$ goes to $+\infty$. Indeed, we define $v_{n}%
(t)=\xi(t-t_{m_{n}})$, which is a solution of (\ref{naci}) with $\beta
_{n}(t)=\beta(t-t_{m_{n}})$. Up to a subsequence, $v_{n}(t)$ converges
uniformly in compact sets to a complete bounded solution $q(t)$ of problem
(\ref{naci}) but replacing $\beta(t)$ by the limit of $\beta_{n}(t)$, denoted
by $\gamma(t)$. It follows that $q(t)\in\alpha_{\xi}(u_{0})\subset
\mathfrak{F}_{j}^{+}\cup\{0\}$, for all $t\in\mathbb{R}$, and that
$q(0)=\widetilde{\varphi}$.

Restricting the system to the subintervals $[0,\frac{\pi}{j}]$, $[\frac{\pi
}{j},\frac{2\pi}{j}]$, $\cdots$, and using comparison, we obtain that $q(t)\in
Y_{j}^{+}\subset\mathcal{O}$, for all $t\geq0$. If $\sigma_{n}$ is bounded,
then we can assume that $\sigma_{n}\rightarrow\sigma$, and then $v_{n}%
(\sigma_{n})\notin\mathcal{O}$ implies that $q(\sigma)\notin\mathcal{O}$,
which is a contradiction.

Further, we define the sequence $u_{n}(t)=\xi(t-t_{m_{n}}+\sigma_{n})$. It is
clear that%
\begin{align*}
u_{n}(t) &  \in\mathcal{O}\text{ for }-\sigma_{n}\leq t<0,\\
u_{n}(0) &  \in\overline{\mathcal{O}},\\
u_{n}(t) &  \not \in \mathcal{O}\text{ for }0\leq t\leq T_{n}.
\end{align*}

This function is a solution of problem (\ref{naci}) with $\overline{\beta}%
_{n}(t)=\beta(t-t_{m_{n}}+\sigma_{n})$.

Up to a subsequence, $u_{n}(t)$ converges uniformly in compact sets to a
complete bounded solution $p(t)$ of problem (\ref{naci}) but replacing
$\beta(t)$ by the limit of $\overline{\beta}_{n}(t)$, denoted by
$\overline{\gamma}(t)$. It follows that $p(t)\in\mathfrak{F}_{j}^{+}\cup \{0\}$ for all
$t$, and that $p(t)\in\overline{\mathcal{O}}$, for all $t\leq0$. Restricting
the system to the subintervals $[0,\frac{\pi}{j}]$, $[\frac{\pi}{j},\frac{2\pi}{j}]$, $\cdots$, we obtain that $p(t)$ is the
unique nondegenerate solution as $t\rightarrow-\infty$, so that $p(t)=\xi
_{j,\bar{\gamma}}^{+}(t)\in Y_{j}^{+}\subset\mathcal{O}$ for all
$t\in\mathbb{R}$.

Furthermore, $u_{n}(0)\rightarrow p(0)$ implies, as $u_{n}(0)\notin%
\mathcal{O}$, that $p(0)\not \in Y_{j}^{+}$, which is a contradiction.

It follows that $\varphi\in Y_{j}^{+}$.

As $\alpha_{\xi}(u_{0})$ is connected, we have that $0\notin\alpha_{\xi}%
(u_{0})$, and then $\alpha_{\xi}(u_{0})\subset Y_{j}^{+}$.

Finally, let us prove (\ref{FormulaAlpha}). Let $\varphi\in\alpha_{\xi}%
(u_{0})\subset Y_{j}^{+}$. We can choose a sequence $t_{n}\rightarrow+\infty$
such that $\xi(-t_{n})\rightarrow\varphi$. Arguing as before we can prove that
$u_{n}(t)=u(t-t_{n})$ converges to a complete bounded solution $p(t)$ of
problem (\ref{naci}) but replacing $\beta(t)$ by the limit of $\beta
_{n}(t)=\beta(t-t_{n})$, denoted by $\gamma(t)\in\mathcal{S}^-(\beta)$.

Since $p(t)\in Y_{j}^{+}$ for all $t\in\mathbb{R}$, it is clear that
$p(t)=\xi_{j,\gamma}^{+}(t)$. Hence, $\varphi=p(0)=\xi_{j,\gamma}^{+}(0)$.$\cqd$

\begin{lemma}
\label{AlfaNonexistence}Let $N^{2}<\lambda<(N+1)^{2}$, $N\geq1$. Then:

\begin{enumerate}
\item There cannot exist a nonzero element $\varphi\in\omega(u_{0})$ such that
$\varphi\in\mathfrak{F}_{j}^{\pm}\cup\{0\}$, with $j\geq N+1.$

\item If $\xi$ is a bounded global solution with initial condition $u_{0},$
there cannot exist a nonzero element $\varphi\in$ $\alpha_{\xi}(u_{0})$ such
that $\varphi\in\mathfrak{F}_{j}^{\pm}\cup\{0\}$, with $j\geq N+1$.
\end{enumerate}
\end{lemma}

\noindent\textbf{Proof:} Suppose by contradiction that such element exists in
$\omega(u_{0})$. Then%
\[
u(t_{n})=T_{\beta}(t_{n},0)u_{0}\rightarrow\varphi\text{ in }C^{1}%
([0,\pi]),\text{ as }n\rightarrow+\infty,
\]
for some sequence $t_{n}\rightarrow+\infty.$

We prove the case where $\varphi\in\mathfrak{F}_{j}^{+}$ with $j$ even. Thus,
$\varphi$ satisfies (\ref{PropFi}).

Let $u_{n}\left(  t\right)  =u\left(  t+t_{n}\right)  $, $t\geq-t_{n}$,
$\beta_{n}\left(  t\right)  =\beta(t+t_{n})$, $t\in\mathbb{R}$. Passing to a
subsequence we obtain that $u_{n}\rightarrow p$, $\beta_{n}\rightarrow\gamma$
uniformly on bounded subsets of $\mathbb{R}$, where $p\left(
\text{\textperiodcentered}\right)  $ is a global bounded solution of the
problem%
\begin{equation}
\left\{
\begin{array}
[c]{c}%
p_{t}=p_{xx}+\lambda p-\gamma(t)p^{3}\text{, }0<x<\pi,\ t\in\mathbb{R},\\
p(t,0)=p(t,\pi)=0,\\
p(0,x)=\varphi(x).
\end{array}
\right. \label{Eqp}%
\end{equation}
From the uniqueness of solutions we have%
\begin{equation}
0=p(t,\frac{\pi}{j})=p(t,\frac{2\pi}{j})=...=p(t,\frac{(j-1)\pi}{j})\text{
}\forall t\geq0.\label{Zeros}%
\end{equation}
Let us prove this fact for $t<0$ as well.

For any $t>0$ take $t_{n_{0}}\geq t$. Then $u(-t)$ is well defined for $n\geq
n_{0}$ and $u_{n}(-t)=u(t_{n}-t)\rightarrow p(-t)$. Therefore, $p(-t)\in
\omega(u_{0})$, so $p(-t)\in\mathfrak{F}_{j}^{+}$ and then (\ref{Zeros}) is
true. The equation (\ref{Eqp}) can be considered separately in each interval
$[0,\frac{\pi}{j}],\ [\frac{\pi}{j},\frac{2\pi}{j}],...$

Consider for instance the first interval. The operator $-\frac{\partial^{2}%
}{\partial x^{2}}$ with Dirichlet boundary conditions in the interval
$[0,\frac{\pi}{j}]$ has the eigenvalues $\lambda_{n}=j^{2}n^{2}$,
$n\in\mathbb{N}$. Since $\lambda\in(N^{2},(N+1)^{2})$ and $j\geq N+1$, every
solution tends to $0$ \cite{Alexandre}. Thus, the unique global bounded
solution to problem (\ref{Eqp}) in $[0,\frac{\pi}{j}]$ is the null solution,
which contradicts the fact that $p\left(  t\right)  $ is a global bounded
nonzero solution. $\cqd$

\begin{lemma}
\label{SymmetryGlobalSol}Let $\lambda>1$. If $\xi($\textperiodcentered$)$ is a
bounded global solution through $u_0$ such that $\omega_{\xi}(u_0),\alpha_{\xi}(u_0)\subset\Xi
_{j}=\{\xi_{j,\gamma}^{\pm}(t):\gamma\in\mathcal{S}^+(\beta)\cup\mathcal{S}^-(\beta),\ t\in\mathbb{R}%
\}$, then for all $a\in\lbrack-\pi,\pi]$ we have:

\begin{enumerate}
\item Either $\rho_{a}\xi(t)-\xi(t)=0$ or $\rho_{a}\xi(t)-\xi(t)\in\Psi$ for any $t\in\mathbb{R}$, where $\xi: \mathbb{R}\rightarrow H^1_P(-\pi,\pi)$ is considered as the odd periodic extension of itself;
\item $\rho_{a}\xi(t)=\xi(t)$ if and only if $\xi_{x}(t,a)=0.$
\end{enumerate}
\end{lemma}

\noindent\textbf{Proof: }Throughout this proof we consider every function
defined in $[-\pi,\pi]$ by taking its odd extension.

Consider two sequences $t_{n},t_{m}\rightarrow+\infty$ such that
\begin{align}
\xi(t_{n}) &  \rightarrow\varphi^{+}\in\Xi_{j},\label{ConvAlfa}\\
\xi(-t_{m}) &  \rightarrow\varphi^{-}\in\Xi_{j}.\nonumber
\end{align}
For any $a\in\lbrack-\pi,\pi]$ the function $w(t)=\rho_{a}\xi(t)-\xi(t)$ is a
global solution of the linear problem%
\[
\left\{
\begin{array}
[c]{c}%
w_{t}=w_{xx}+r(t,x)w,\ x\in (-\pi,\pi),\ t\in\mathbb{R},\\
w(-\pi,t)=w(\pi,t), \ w_x(-\pi,t)=w_x(\pi,t), \ t\in\mathbb{R},\\
w(0,x)=(\rho_{a}\xi(0)-\xi(0))(x),\ x\in (-\pi,\pi),
\end{array}
\right.
\]
where $r(t,x)= \lambda - \beta(t)\left[\frac{(\rho_a\xi)^3-\xi^3}{\rho_a\xi-\xi}\right]$.
For simplicity of notation we omit the dependence of $w$ on $a$.
First, take $a$ such that $\rho_{a}\varphi^{+}-\varphi^{+}\not =0$ and
$\rho_{a}\varphi^{-}-\varphi^{-}\not =0$. In such a case the lap numbers of
these functions satisfy:%
\[
\ell(\rho_{a}\varphi^{+}-\varphi^{+})=\ell(\rho_{a}\varphi^{-}-\varphi^{-}%
)<\infty\text{.}%
\]
Indeed, for any $\varphi\!\in\!\mathfrak{F}_{j}(-\frac{\pi}{2j})\,\cup\,
\mathfrak{F}_{j}(\frac{\pi}{2j})$ such that $f_{a}=\rho_{a}\varphi
-\varphi\not =0$ and $\rho_{a}\varphi\not \in \mathfrak{F}_{j}(-\frac{\pi
}{2j})\,\cup\,\mathfrak{F}_{j}(\frac{\pi}{2j})$ the zeros of the function
$f_{a}$ are located at the points%
\[
-\pi+a,\ -\frac{(j-1)\pi}{j}+a,\ ...,\ -\frac{\pi}{j}+a,\ a,\ \frac{\pi}%
{j}+a,...,\frac{(j-1)\pi}{j}+a,
\]
that is, there are $2j$ zeros. If $f_{a}\not =0$ and $\rho_{a}\varphi
\in\mathfrak{F}_{j}(-\frac{\pi}{2j})\cup\mathfrak{F}_{j}(\frac{\pi
}{2j})$, then the zeros are located at%
\[
-\pi,-\frac{(j-1)\pi}{j},...,-\frac{\pi}{j},0,\frac{\pi}{j},...,\frac
{(j-1)\pi}{j},\pi,
\]
that is, there are $2j+1$ zeros.

Since $w(t_{n})\rightarrow\rho_{a}\varphi^{+}-\varphi^{+},$ $w(-t_{m}%
)\rightarrow\rho_{a}\varphi^{-}-\varphi^{-}$, we obtain that%
\begin{align*}
\ell(w(t)) &  =\ell(\rho_{a}\varphi^{+}-\varphi^{+})\text{ }\forall t\geq
t_{_{1}},\\
\ell(w(t)) &  =\ell(\rho_{a}\varphi^{-}-\varphi^{-})\text{ }\forall
t\leq-t_{_{2}},
\end{align*}
Hence, choosing an arbitrary $t^{\ast}\geq\max\{t_{_{1}},t_{_{2}}\}$ we obtain%
\[
\ell(w(t^{\ast}))=\ell(w(-t^{\ast}))\text{.}%
\]
Using again point (iii) of Lemma \ref{angenent} it follows that $w(t)\in\Psi$
for all $t\in\mathbb{R}.$

Consider now $a$ such that $\rho_{a}\varphi^{+}-\varphi^{+}=\rho_{a}%
\varphi^{-}-\varphi^{-}=0$. Then%
\[
w\left(  t\right)  \rightarrow0\text{ as }t\rightarrow\pm\infty\text{.}%
\]
But we proved in Section \ref{S3}, before of the Theorem \ref{lambdaentre1e4}, that there cannot exist homoclinic solutions. Then $w(t)\equiv 0$ and
the first statement is proved. The second one is an easy consequence of the first one.$\cqd$

\begin{lemma}
\label{HomoclinicNonexistence}Let $\lambda>1$. If $\xi($\textperiodcentered$)$
is a bounded global solution through $u_0$ for $\beta\left(  \text{\textperiodcentered
}\right)$ such that $\omega_{\xi}(u_0),\alpha_{\xi}(u_0)\subset\Xi_{j}=\{\xi
_{j,\gamma}^{\pm}(t):\gamma\in\mathcal{S}^+(\beta)\cup\mathcal{S}^-(\beta),\ t\in\mathbb{R}\}$, $j\geq
1$, then either $\xi=\xi_{j,\beta}^{+}$ or $\xi=\xi_{j,\beta}^{-}.$
\end{lemma}

\noindent\textbf{Proof:} Let $t\in\mathbb{R}$ be arbitrary. We consider the odd
extension of $\xi\left(  t\right)  $ over $[-\pi,\pi]$ and choose $x^{\ast}%
\in\lbrack-\pi,\pi]$ such that $\xi_{x}(t,x^{\ast})\not =0$. For example, let
$\xi_{x}(t,x^{\ast})>0$. Let $I=\left(  x_{0},x_{1}\right)  $ be the maximal
interval containing $x^{\ast}$ where $\xi_{x}(t,x)$ for any $x\in I$. $\xi
(t,$\textperiodcentered$)\in C^{1}([-\pi,\pi])$ implies that $\xi_{x}%
(t,x_{0})=\xi_{x}(t,x_{1})=0$. Then using Lemma \ref{SymmetryGlobalSol} we get
$\rho_{x_{0}}\xi(t)=\rho_{x_{1}}\xi(t)=\xi(t)$. Therefore, there exists
$m\in\mathbb{N}$ such that $x_{1}-x_{0}=\frac{\pi}{m}$, which gives that either
$\xi(t)\in\mathfrak{F}_{m}(\frac{\pi}{2m})$ or $\xi(t)\in\mathfrak{F}%
_{m}(-\frac{\pi}{2m})$.

Consider two sequences $t_{n},t_{m}\rightarrow+\infty$ such that
\begin{align}
\xi(t_{n}) &  \rightarrow\varphi^{+}\in\Xi_{j},\label{ConvAlfa2}\\
\xi(-t_{m}) &  \rightarrow\varphi^{-}\in\Xi_{j}.\nonumber
\end{align}
From (\ref{ConvAlfa2}) we deduce that%
\begin{align*}
\ell(\xi(t)) &  =\ell(\varphi^{+})\text{ }\forall t\geq t_{_{1}},\\
\ell(\xi(t)) &  =\ell(\varphi^{-})\text{ }\forall t\leq-t_{_{2}},
\end{align*}
Hence, choosing an arbitrary $t^{\ast}\geq\max\{t_{_{1}},t_{_{2}}\}$ we obtain%
\[
\ell(\xi(t^{\ast}))=\ell(\xi(-t^{\ast}))=2j+1\text{.}%
\]
Therefore, since the number of zeros of $\xi\left(  t\right)  $ is
non-increasing, $m=j$ for any $t\in\mathbb{R}$. Also, we observe that as a
consequence of Lemma \ref{angenent} the function $\xi(t)$ cannot jump from $\mathfrak{F}%
_{j}(\frac{\pi}{2j})$ into $\mathfrak{F}_{j}(-\frac{\pi}{2j})$ or
viceversa, that is, either $\xi(t)\in\mathfrak{F}_{j}(\frac{\pi}{2j})$,
for all $t\in\mathbb{R}$, or $\xi(t)\in\mathfrak{F}_{j}(-\frac{\pi}{2j})$,
for all $t\in\mathbb{R}$.

Now we can consider the solution $\xi(t)$ separately in each subinterval
$[0,\frac{\pi}{j}],\ [\frac{\pi}{j},\frac{2\pi}{j}],...$ Since $\alpha_{\xi
}(u_0)\subset\Xi_{j}$, $\xi(t)$ restricted to the interval $[0,\frac{\pi}{j}]$ is
non-degenerate at $-\infty$. Moreover, by Lemma \ref{AlfaNonexistence} we know
that $\lambda>j^{2}$. As the eigenvalues of the operator $-\partial_{x}^{2}$
on $[0,\frac{\pi}{j}]$ are $\lambda_{n,\frac{\pi}{j}}=n^{2}j^{2}$, we obtain
that $\lambda>\lambda_{1,\frac{\pi}{j}}$. If for example $\xi(t)\in
\mathfrak{F}_{j}(-\frac{\pi}{2j})$, this implies that $\xi\left(
t\right)  \mid_{\lbrack0,\frac{\pi}{j}]}=\xi_{1,\beta,\frac{\pi}{j}}^{+}$,
that is, it coincides with the unique positive non-degenerate solution of
problem (\ref{naci}) on the interval $[0,\frac{\pi}{j}]$. Repeating the same
argument in each interval we finally prove the equality $\xi=\xi_{j,\beta}%
^{+}$. If $\xi(t)\in\mathfrak{F}_{j}(\frac{\pi}{2j})$, in the same way we
have $\xi=\xi_{j,\beta}^{-}.$ $\cqd$

\section{Connections between the zero and the ``non-autonomous equilibria"}
\label{S4}

In this section, for $N^2<\lambda<(N+1)^2$, $N \geq 1$, we will prove the existence of global bounded solutions, $\zeta_j^{\pm}(t)$, which connect the trivial solution $\xi_0 \equiv 0$ to the non-autonomous equilibria  $\xi_j^{\pm}(t)$, $1\leq j \leq N$. We remember that an analogous fact also occurs in the autonomous case. This makes us believe that the solutions $\xi_j^\pm(t)$ play, for the non-autonomous problem, the same role that the equilibria $\phi_{j,\beta_i}^\pm$ in the autonomous case.

Resuming our goal, we observe that, due to the symmetries of the equation, it is enough to verify the existence of $\zeta_j^+(t)$, so the solution $\zeta_j^-(t) := - \zeta_j^+(t)$ will connect $\xi_0$ to $\xi_j^-(t) = - \xi_j^+(t)$.

Further, we will see that it is sufficient to show the existence of a positive solution $\zeta_1^+(t)$, $t \in\mathbb{R}$, in the unstable manifold of zero, that the existence of the others solutions $\zeta_j^+(t)$, $2 \leq j \leq N$, will be a consequence. The global solutions $\zeta_j^+(t)$ will get out from zero and come in the strip $Y_j^+$, for $t$ sufficiently large. It follows from the Theorem \ref{Rodrigues-Bernal} that $\zeta_j^+(t)- \xi_j^+(t) \to 0$, as $t \to +\infty$.


\begin{theorem}\label{Theo-from-zero}
If $N^2< \lambda <(N+1)^2$, for each $1\leq j\leq N$ there are global solutions $\zeta_j^\pm$ of \eqref{naci} such that  $0\stackrel{t\to -\infty}{\longleftarrow} \zeta_j^\pm (t)\stackrel{t\to +\infty}{\longrightarrow}\xi_j^\pm$ where $\xi_j^\pm$ are the non-autonomous equilibria.
\end{theorem}

The proof of this theorem is a consequence of the symmetry properties of solutions and of the following proposition.

\begin{proposition}
\label{lambda<4}
If $1< \lambda \neq n^2$, $n\geq 2$, then the solution $\zeta (t):= \zeta_1^+(t)$, given by Theorem \ref{lambdaentre1e4}, is positive for any $t \in \mathbb{R}$.
\end{proposition}
 \noindent{\bf Proof:} If $1< \lambda <4$, consider the space $X_1$ generated by the first eigenfunction $\sin(x)$. The local unstable manifold of zero can be obtained as the graph of a map $\sigma : \mathbb{R}\times X_1 \rightarrow H_0^1(0,\pi)$, i.e,
$$W_{loc}^u(0) \subset \{(t, u_1, u_2): u_2= \sigma(t, u_1), t\in \mathbb{R}, u_1 \in X_1\}.$$
Furthermore, the graph of $\sigma(t,\cdot)$ is tangent to the subspace $X_1$ at the origin. Equivalently, given $\tau \in \mathbb{R}$ and $s_0>0$ sufficiently small, we have that the global solution which constitutes the unstable manifold of the zero can be written as $\zeta(t)= \zeta_1(t)+ \sigma(t, \zeta_1(t))= \zeta_1(t) + \zeta_2(t)$, $t\leq \tau$, where $\zeta(\tau)= s_0\sin(x) + \sigma(\tau, s_0\sin(x))$ and
$$\frac{\|\zeta_2(t)\|_{H_0^1}}{\|\zeta_1(t)\|_{H_0^1}} = \frac{\|\sigma(t, \zeta_1(t))\|_{H_0^1}}{\|\zeta_1(t)\|_{H_0^1}} \rightarrow \;0,\;\;\mbox{as}\;t \rightarrow - \infty.$$
Hence
$$\zeta(t,x)= s_1(t)\sin(x)+ \sum_{n\geq 2}s_n(t)\sin(nx)$$
where $s_1(t)> 0$ and  $\zeta(t,x)>0$ for all $0<x<\pi$ and $t \in \mathbb{R}$, as  consequence of Lemma \ref{angenent}. For the general case the proof is similar with the only difference that there is a local invariant manifold, given as a graph over $X_1$ and tangent to it at the origin, within the higher dimensional unstable manifold. $\cqd$

\bigskip

The proof of Theorem \ref{Theo-from-zero} now follows from symmetry properties of solutions as sketched in the case $4 < \lambda \leq 9$. In this case, $N=2$ and there are four non-autonomous equilibria $\xi_1^\pm(t)$ and $\xi_2^\pm(t)$, $t \in \mathbb{R}$. Remember that, for any $t \in \mathbb{R}$,
 $$\begin{array}{l}
    \xi_2^+(t, \frac{\pi}{2})=0\\
  \xi_2^+(t,x)>0,\;\;0<x<\frac{\pi}{2} \\
   \xi_2^+(t,\pi -x)= - \xi_2^+(t,x)= \xi_2^-(t,x),\;\;0\leq x \leq \pi
 \end{array}
 $$
  So, $\xi_2^+(t)$ restricted to the interval $[0,\frac{\pi}{2}]$ is a global solution of the non-autonomous equation, with  homogeneous Dirichlet conditions on the interval $[0,\frac{\pi}{2}]$, and the operator $-\partial_x^2$, under these conditions, has eigenvalues $\lambda_{n,\frac{\pi}{2}}= 4n^2$, $n=1,2,\cdots$. Thus, $4 = \lambda_{1,\frac{\pi}{2}} < \lambda \leq 9 < \lambda_{2,\frac{\pi}{2}}$, this is, $\lambda$ is  between the first and second eigenvalue of the operator $-\partial_x^2$ on the interval $[0,\pi/2]$ and, therefore, we can  apply the Preposition \ref{lambda<4} 
to ensure the existence of a positive global bounded solution on the interval $(0,\frac{\pi}{2})$, which we will call $\zeta_{1,\frac{\pi}{2}}^+(t)$ and such that   $$\zeta_{1,\frac{\pi}{2}}^+(-t) \rightarrow 0\;\;\; \mbox{and}\;\;\; \zeta_{1,\frac{\pi}{2}}^+(t)- \xi_{2}^+(t)\Big|_{[0,\frac{\pi}{2}]} \rightarrow 0, \quad \hbox{as}\quad t\to \infty.$$
Note that, when we have the autonomous problem on the interval $[0,\frac{\pi}{2}]$, the equilibrium solutions are given by
  $$\phi_{k,\beta ,\frac{\pi}{2}}^\pm = \phi_{2k,\beta}^\pm\Big|_{[0,\frac{\pi}{2}]}.$$
  Consequently, the non-autonomous equilibria of the non-autonomous equation on the interval $[0,\frac{\pi}{2}]$ are given by $\xi_{k,\frac{\pi}{2}}^{\pm}(t)= \xi_{2k}^{\pm}(t)\Big|_{[0,\frac{\pi}{2}]}$.
\noindent
  Thus, the function defined on $\mathbb{R}\times [0,\pi] $ by
  $$\zeta_2^+(t,x)=
  \left\{
  \begin{array}{l}
    \zeta_{1,\frac{\pi}{2}}^+(t,x), \;0\leq x \leq \frac{\pi}{2}\\
     - \zeta_{1,\frac{\pi}{2}}^+(t, \pi -x), \;\frac{\pi}{2}< x \leq \pi
  \end{array}
  \right.
  $$
   is a solution of equation (\ref{naci}) which goes out of zero and tends to the solution $\xi_2^+(t)$. Also, the zeroes of $\zeta_2^+(t)$ are fixed: $0$, $\pi/2$ and $\pi$. $\cqd$

\section{Final comments and further problems}\label{comments}

As a consequence of the results in this paper, we have completely characterised the skew product attractor for \eqref{naci}. If $N^2< \lambda <(N+1)^2$, $\Sigma=cl_{C(\mathbb{R},\mathbb{R})}\{\beta(\cdot+t):=\theta_t\beta: t\in \R^+\}$, $\mathcal{S^{+}(\beta)}$ (with respect to the metric of the uniform convergence in bounded subsets or $\R^+$) is the global attractor of $\{\theta_t: t\in \R^+$ in $\Sigma$, and defining ${\mathbf{\Xi}}_j^\pm = \{(\xi_{j,\gamma}^{\pm}(t),\gamma):\gamma \in \mathcal{S}^{+}(\beta),\; t \in \mathbb{R}\}$, $1\leq j  \leq N$, ${\mathbf{\Xi}}_{N+1}^+=\{0\}$, then $\{{\mathbb{Z}}_1, \cdots, {\mathbb{Z}}_{2N+1}\}$, where $\mathbb{Z}_{2j-1}={\mathbf{\Xi}}_j^+$, $1\leq j\leq N+1$, $\mathbb{Z}_{2j}={\mathbf{\Xi}}_j^-$, $1\leq j\leq N$,  is a Morse decomposition for the skew product semiflow $\Pi(t):H^1_0(0,\pi)\times \Sigma \to H^1_0(0,\pi)\times \Sigma$ given by $\Pi(t)(u_0,\gamma)=(T_\gamma(t,0)u_0,\theta_{t}\gamma)$ and $T_\gamma(t,0)u_0$ is the solution at time $t$ of \eqref{naci} with $\beta$ replaced by $\gamma$ that have started at $u_0$. In fact Theorem \ref{charac-omega} and Theorem \ref{charac-alpha} ensure that all global solutions either connect two of the invariant sets $\mathbb{Z}_{j}$ or are contained in one of them, these sets are invariant, closed and disjoint and Lemma \ref{HomoclinicNonexistence} ensures the non-existence of homoclinic structures, so they are maximal invariant.

This give a non-trivial example of skew-product semigroup with gradient structure. We end this paper by expressing our belief that there are yet some needed work for this very nice example.

\begin{conjecture}
Concerning \eqref{naci} and for $N^2< \lambda <(N+1)^2$, we conjecture that,
\begin{enumerate}
\item for each $j= 1, 2,\cdots ,N$, the non-autonomous equilibria $\xi_j^\pm(t)$ are connected  to all the equilibria $\xi_k^\pm(t)$, $k=1,\cdots,j-1$,
\item that the pullback attractor $\mathcal{A}(t)=\overline{W^u(0)(t)}=\bigcup_{j=1}^{N} W^u(\xi_j^\pm)(t) \cup W^u(0)(t)$,
\item that the non-autonomous equilibria $\xi_j^\pm$, $1\leq j\leq N$, are hyperbolic and
\item the unstable and stable manifolds of the non-autonomous equilibria intersect transversally along a connection.
\end{enumerate}
All of these are true for the autonomous case $\beta(t)\equiv \hbox{const}$.
\end{conjecture}

We also comment that the exact form of $f$ was never used. In fact we only need that $f$ be an odd smooth function.

\bigskip

%


\begin{thebibliography}{99}
\bibitem{Angenent} S. Angenent, The zero set of a solution of a parabolic problem, \textit{Journal f\"{u}r die reine unde angewandte Mathematik}, {\bf 390} (1988), 79-96.

\bibitem{Angenent2} S. Angenent,
The Morse-Smale property for a semilinear parabolic equation.
\textit{J. Differential Equations} \textbf{62} (1986), no. 3, 427-442.

\bibitem{ACCL} E. R. Arag\~{a}o-Costa, T. Caraballo, A. N. Carvalho and  J. A. Langa, Stability of gradient semigroups under perturbations, {\sl Nonlinearity} {\bf 24}  2099-2117 (2011).


\bibitem{BV} A. V. Babin and M. I. Vishik, \textit{Attractors in Evolutionary Equations} Studies in Mathematics and its Applications \textbf{25},
North-Holland Publishing Co., Amsterdam, 1992.

\bibitem{BCLR} M. C. Bortolan, A. N. Carvalho, J. A. Langa and G. Raugel, Non-autonomous perturbations of Morse-Smale semigroups: stability of the phase diagram,  Preprint.

\bibitem{CLR12-1} A. N. Carvalho, J. Langa and J. C. Robinson, Structure and bifurcation of pullback attractors in a non-autonomous Chafee-Infante equation, \textit{Proceedings of the American Mathematical Society} {\bf 140} 2357-2373 (2012).

\bibitem{CLR12} A. N. Carvalho, J. Langa and J. C. Robinson, Attractors for infinite-dimensional non-autonomous dynamical systems, Applied Mathematical Sciences {\bf 182}  Springer-Verlag 2012.

\bibitem{CL07} A. N. Carvalho and J. A. Langa, Non-autonomous perturbation of autonomous semilinear differential equations: Continuity of local stable and unstable manifolds, \textit{J. Differential Equations} {\bf 233} 622-653 (2007).

\bibitem{CLRS} A. N. Carvalho, J. A. Langa, J. C. Robinson and A. Su\'{a}rez, Characterization of non-autonomous
attractors of a perturbed infinite-dimensional gradient system,
\textit{J. Differential Equations}  \textbf{236}  (2007),  no. 2, 570--603.

\bibitem{CL09} A. N. Carvalho and J. A. Langa, An extension of the concept of
gradient semigroups which is stable under perturbation, \textit{J.
Differential Equations} \textbf{246} (2009), 2646--2668.

\bibitem{Alexandre}{\sc A.N.Carvalho, J. A. Langa and J. C. Robinson}, {\it Structure and bifurcation of pullback attractors in a non-autonomous Chafee-Infante equation} {\it Proceedings of the American Mathematical Society}, {\bf 140} (2012), 2357-2373.

\bibitem{chen}{\sc  X.Y. Chen and H. Matano}, {\it Convergence, asymptotic periodicity and finite-point blow-up in one dimensional semilinear heat equations}, JDE {\bf 78} (1989), 160-190 .

\bibitem{CI} N.~Chafee and E.~F. Infante. A bifurcation problem for a nonlinear partial differential equation
  of parabolic type. {\em Applicable Anal.}, 4:17--37, 1974/75.

\bibitem{hale}{\sc  M. Chen, X-Y. Chen and J. K. Hale}, {\it Strutural stability for time-periodic one-dimensional parabolic equations}, JDE {\bf 96} (1992), 335-418 .

\bibitem{FRJDE1996} B. Fiedler and C. Rocha, Heteroclinic orbits of semilinear parabolic equations, \textit{J. Differential Equations} \textbf{125} (1996), no. 1, 239-281.

\bibitem{Hale}  J. Hale, \textit{Asymptotic Behavior of Dissipative Systems}
(Providence: Math. Surveys and Monographs, A.M.S. 1998).

\bibitem{Henry}  D. Henry, \textit{Geometric theory of semilinear parabolic
equations, }Lecture Notes in Mathematics 840, Berlin: Springer,
1981.

\bibitem{Henry2}  D. Henry, Some infinite-dimensional Morse-Smale systems defined by parabolic partial differential equations. \textit{J. Differential Equations} \textbf{59} (1985), no. 2, 165-205.


\bibitem{LS} J. A. Langa and A. Suarez, Pullback permanence for non-autonomous partial differential
equations, Electronic J. Differential Equations 2002, No. 72.

\bibitem{Matano} H. Matano,  Nonincrease of the lap-number of a solution for a one-dimensional semilinear parabolic equation,  \textit{J. Fac. Sci. Univ. Tokyo Sect. IA Math.}  29  (1982), no. 2, 401--441.

\bibitem{RB-VL} A. Rodriguez-Bernal and A. Vidal-Lopez, Existence, uniqueness and attractivity properties of positive complete trajectories for non-autonomous reaction-diffusion problems. {\it Dis. Cont. Dyn. Sys.} Series A,{\bf 18} (2007), 537-567.

\end{thebibliography}
\end{document}